\documentclass[english,10pt]{amsart}

\usepackage{enumerate}
\usepackage{amssymb}
\usepackage{array}
\usepackage{graphicx}
\usepackage{multirow}
\usepackage{mathrsfs}
\usepackage{multicol}
\usepackage[a4paper,margin=3cm]{geometry}
\usepackage{mathtools}
\usepackage{float}
\usepackage{color}
\usepackage{longtable}
\usepackage[dvipsnames,svgnames,table]{xcolor}
\usepackage{comment}
\usepackage{pgf,tikz,pgfplots}
\usepackage{tikz-cd}
\usetikzlibrary{arrows}
\pgfplotsset{compat=1.15}

\usepackage[curve,all]{xy}
\xyoption{matrix}
 \xyoption{curve}
 \xyoption{color}
 \xyoption{line}
\xyoption{arc}

\usepackage[shortlabels]{enumitem}
\setlist[enumerate]{leftmargin=7mm,topsep=0pt,itemsep=-1ex,partopsep=1ex,parsep=1ex,label=\rm{(\roman*)}}
\setlist[itemize]{leftmargin=5mm,topsep=0pt,itemsep=-1ex,partopsep=1ex,parsep=1ex,label=\raisebox{0.25ex}{\tiny$\bullet$}}
\xdefinecolor{tu-blue}{RGB}{19,83,138}
\usepackage[backref]{hyperref}

\usepackage[nameinlink,capitalize,noabbrev]{cleveref}

\theoremstyle{plain}
\newtheorem{theorem}{Theorem}[section]

\newtheorem*{theoremaux}{Theorem \theoremauxnum}
\gdef\theoremauxnum{1}

\newenvironment{theoremff}[2][]{
  \def\theoremauxnum{\ref{#2}}
  \begin{theoremaux}[#1]
}{%
  \end{theoremaux}
}
\newtheorem*{main-theorem}{Main Theorem}
\newtheorem{proposition}[theorem]{Proposition}
\newtheorem*{propositionaux}{Proposition \propositionauxnum}
\gdef\propositionauxnum{1}

\newtheorem{lemma}[theorem]{Lemma}
\newtheorem*{lemmaaux}{Lemma \lemmaauxnum}
\gdef\lemmaauxnum{1}

\newtheorem{corollary}[theorem]{Corollary}

\newtheorem*{key-problem}{Key Problem}
\newtheorem*{theorem-intro}{Theorem}
\newtheorem*{proposition-intro}{Proposition}
\newtheorem*{corollary-intro}{Corollary}

\theoremstyle{definition}
\newtheorem{definition}[theorem]{Definition}

\newtheorem{example}[theorem]{Example}

\theoremstyle{remark}
\newtheorem{remark}[theorem]{Remark}

\makeatletter
\@namedef{subjclassname@2020}{%
  \textup{2020} Mathematics Subject Classification}
\makeatother
\newcommand{\incl}[1][r]{\ar@<-0.2pc>@{^(-}[#1] \ar@<+0.2pc>@{-}[#1]}

\newcommand{\hs}{\kern 0.8pt}

\renewcommand{\H}{{\mathrm{H}}}

\newcommand{\Loc}{\mathrm{Loc}}

\newcommand{\Hom}{\mathrm{Hom}}

\renewcommand{\div}{\mathrm{div}}
\newcommand{\E}{\mathfrak{E}}

\renewcommand{\H}{\mathrm{H}}
\newcommand{\Supp}{\mathrm{Supp}}
\newcommand{\Spec}{\mathrm{Spec}}

\newcommand{\pr}{\mathrm{pr}}
\newcommand{\Q}{\mathbb{Q}}

\renewcommand{\H}{\mathrm{H}}

\newcommand{\A}{\mathbb{A}}

\newcommand{\D}{\mathfrak{D}}
\renewcommand{\S}{\mathfrak{S}}

\newcommand{\N}{\mathbb{N}}

\newcommand{\Z}{\mathbb{Z}}

\newcommand{\Gm}{\mathbb{G}_\mathrm{m}}

\newcommand{\Div}{\mathrm{Div}}

\newcommand{\aff}[1]{#1_{\mathrm{aff}}}

\DeclareMathOperator{\rk}{rk}

\DeclareMathOperator{\Sing}{Sing}

\title[]{Equivariant resolution of singularities for $T$-varieties in positive characteristic}
\date{}

\author{Gary Mart\'inez-N\'unez}
\address{Univ. Lille, CNRS, UMR 8524 - Laboratoire Paul Painlev\'e, F-59000 Lille, France}
\email{gary-antonio.martinez-nunez@univ-lille.fr}

\begin{document}

\begin{abstract}
Let the ground field be algebraically closed of positive characteristic. We prove that every variety endowed with an algebraic torus action whose general orbits have codimension at most three admits an equivariant resolution of singularities. The proof combines the description of such varieties in terms of proper polyhedral divisors and divisorial fans, the existence of log resolutions in dimension at most three, and the resolution of toroidal embeddings established by Kempf, Knudsen, Mumford and Saint-Donat.
\end{abstract}
\maketitle


\section{Introduction}
A \emph{variety} is a separated geometrically integral scheme of finite type over a field. A \emph{resolution of singularities} (or simply a \emph{desingularization}) of a variety $X$ is a proper birational morphism $\Psi:W\to X$ from a smooth variety $W$. When $X$ is a $G$-variety, for $G$ an algebraic group, a resolution of singularities $\Psi:W\to X$ is called \emph{equivariant} if $W$ is also a $G$-variety and $\Psi$ is equivariant.

In 1964, Hironaka \cite{Hir64} proved that every variety over a field of characteristic zero admits a resolution of singularities. In positive characteristic, the general problem remains open, although curves, surfaces, and threefolds admit resolutions of singularities. The surface case was established by Abhyankar \cite{Abh56} and Lipman \cite{Lip78}. Abhyankar \cite{Abh66} proved resolution of singularities for threefolds over algebraically closed fields of characteristic $p> 5$ (see Cutkosky \cite{Cut09} for a concise proof), and Cossart and Piltant \cite{CP08,CP09} removed the restriction on the characteristic; they later extended their result to arbitrary quasi-excellent schemes of dimension at most three, together with a log-resolution statement \cite[Theorem 1.1]{CP19}. 
For a survey of the state of the art in resolution of singularities, we refer the reader to Spivakovsky \cite{Spi20}.

Equivariant resolution of singularities exists in characteristic zero (see \cite[Proposition 3.9.1 and Theorem 3.36]{Kol07}). In positive characteristic, the general problem remains open, but several important classes of varieties are known to admit equivariant resolutions. These include toric varieties \cite[Theorem 11, p. 34]{KKMSD73}, \cite[Section 2.6]{Ful93}, as well as the embedded toric and binomial varieties treated by Bierstone and Milman \cite{BM06}. Brion \cite[Proposition 2.4.4(iii)]{Bri26} shows that every $G$-variety such that every non-empty closed $G$-stable subscheme has codimension at most two admits a strong $G$-desingularization (see \cite[Definition 2.4.1]{Bri26}). 

\subsection*{Main results.}
Throughout the rest of the paper, $k$ is an algebraically closed field of positive characteristic and $T$ is an algebraic torus over $k$. The complexity of a $T$-variety is the codimension of the general orbit. Our main result claims the following.

\begin{theorem}\label{theorem: main theorem affine}
Every $T$-variety of complexity at most three admits an equivariant resolution of singularities.
\end{theorem}

As a consequence of \cref{theorem: main theorem affine}, every four-dimensional $T$-variety admits an equivariant resolution of singularities. The restriction on the complexity arises from the lack of general results on resolution of singularities in arbitrary dimensions. It is worth emphasizing that, once a general resolution of singularities theorem is established, the existence of equivariant resolution in arbitrary dimension would follow automatically from our work.

Since $T$ is a smooth algebraic group, every $T$-variety admits an equivariant normalization. Therefore, the problem can be reduced to the desingularization of normal $T$-varieties. These varieties were studied by Altmann and Hausen \cite{AH06} and by S\"u\ss\ \cite{AHS08}, who provided and algebro-geometric combinatorial description over algebraically closed fields of characteristic zero. These results were subsequently extended in \cite{MN25a} and \cite{MN25b} to arbitrary fields and in any characteristic. Such descriptions have being used in \cite{LS13} to study the singularities of $T$-varieties and their resolutions over algebraically closed fields of characteristic zero. In what follows, we use them to construct equivariant resolution of singularities by means of \emph{toroidal embeddings} \cite[Definition 1, p.~54]{KKMSD73}.

\subsection*{Sketch of the proof.} As said above, we reduced ourselves to the case of normal $T$-varieties. By \cite[Proposition 4.17]{MN25a}, for any affine normal $T$-variety $X$, there exists a normal $T$-variety $\tilde{X}$ together with a proper birational $T$-equivariant morphism $r:\tilde{X}\to X$. Such a normal $T$-variety $\tilde{X}$ is not unique. In \cref{Section: partial desingularization non-affine} we prove that an analogous construction is also possible for non-affine normal $T$-varieties. The variety $\tilde{X}$ is determined by the corresponding algebro-geometric combinatorial data: \emph{proper polyhedral divisors} in the affine case and \emph{divisorial fans} otherwise.

We then introduce a particular class of these objects, namely, \emph{simple normal crossing} proper polyhedral divisors (resp. divisorial fans), and prove that they can be constructed in our setting. When $\tilde{X}$ arises from one of these objects, it admits a $T$-stable open subvariety $U$ such that $U\subset \tilde{X}$ is a \emph{toroidal embedding}. In \cref{Section: toroidal embedding and equivariant resolution}, we prove that every toroidal embedding $U\subset \tilde{X}$, with $U$ being $T$-stable, admits an equivariant desingularization. This establishes \cref{theorem: main theorem affine}.

\subsection*{Notation and conventions.}
As above, $k$ is an algebraically closed field of positive characteristic. A \emph{variety} over $k$ is a separated integral scheme of finite type over $k$. Since $k$ is perfect, a variety is regular if and only if it is smooth. We denote by $N$
the cocharacter lattice and by $M:=\Hom_{\Z}(N,\Z)$ the character lattice of $T$ and $N_{\Q}$ the associated $\Q$-vector spaces. By a $T$-variety, we mean a variety endowed with an effective action of $T$. Toric varieties are always assumed to be normal.

\subsection*{Layout of the paper.}
In \cref{Section: recall T-varieties}, we briefly recall the theory of affine normal $T$-varieties. We begin with basic notions from convex geometry, then recall \emph{proper polyhedral divisors} and their associated varieties, and finally discuss torus invariant divisors.

In \cref{Section: equivariant resolution}, we recall the notion of a \emph{toroidal embedding}, prove the equivariance of canonical resolutions of $T$-stable toroidal embeddings, introduce \emph{simple normal crossing pp-divisors}, prove the local structure theorem, study pull-backs of pp-divisors and log resolutions of the base, and prove \cref{theorem: main theorem affine}.

In \cref{Section: resolution general case}, we first recall the theory of non-affine normal $T$-varietie and their description in terms of \emph{divisorial fans}. Then, we construct another normal $T$-variety $\tilde{X}$ with a proper birational $T$-equivariant morphism $r:\tilde{X}\to X$, introduce simple normal crossing divisorial fans, and conclude the section with a proof of \cref{theorem: main theorem affine} in the general case.

Finally, in \cref{Section: combinatorial resolution}, we give a combinatorial description of the resolution and conclude with an example and some open questions.

\subsection*{Acknowledgements.}
The author acknowledges the CDP C2EMPI, the University of Lille, the Initiative of Excellence of the University of Lille, and the European Metropolis of Lille for their funding and support of the R-CDP-24-004-C2EMPI project. The author also thanks \'Alvaro Liendo for many enlightening conversations, and Ronan Terpereau and Patrick Popescu-Pampu for their helpful suggestions regarding the introduction.

\subsection*{AI disclosure.} The mathematical ideas, results, and proofs in this paper are the work of the author and arose as a natural continuation of \cite{LS13}. During the preparation of the manuscript, the author used OpenAI's GPT-6.1 Terra as a tool for mathematical discussions and editorial assistance. This tool was used to examine proofs and hypotheses, explore simplifications, identify minor errors, and improve the organization and exposition of the text. All suggestions were critically assessed and independently verified by the author, who takes full responsibility for the content of the paper.

\section{Brief recall on normal \texorpdfstring{$T$}{T}-varieties}\label{Section: recall T-varieties}

This section summarizes \cite{AH06} and \cite{MN25a}, together with some elementary facts that will be used later.

\subsection{Semi-projective varieties.} The \emph{affinisation} of a variety $Y$ over $k$ is the scheme $\aff{Y}:=\Spec(\H^{0}(Y,\mathscr{O}_{Y}))$, and the morphism $r_{Y}:Y\to \aff{Y}$ corresponding to the identity in $\Hom_{\mathrm{Ring}}(\H^{0}(Y,\mathscr{O}_{Y}),\H^{0}(Y,\mathscr{O}_{Y}))$ is called the \emph{affinisation morphism}. A variety $Y$ over $k$ is called \emph{semi-projective} if $\aff{Y}$ is of finite type and $r_{Y}$ is projective.

\subsection{Convex geometry.} Recall that $N$ is a lattice and $M:=\Hom_{\Z}(N,\Z)$ is its dual. By a \emph{cone} $\omega\subset N_{\Q}$ we mean a convex rational polyhedral cone, which is called \emph{pointed} if it contains no lines. We denote by $\omega^{\vee}\subset M_{\Q}$ its \emph{dual cone}, defined by the perfect pairing $\langle,\rangle:M\times N\to \Z$, which extends linearly to the respective $\Q$-vector spaces, and by $\omega(1)$ its set of rays. For references, see \cite{Ful93} or \cite{CLS11}.

A \emph{polytope} $\Pi\subset N_{\Q}$ is the convex hull of finitely many points in $N_{\Q}$. By a \emph{polyhedron} $\Delta\subset N_{\Q}$ we mean a convex subset of $N_{\Q}$ that can be written as $\Delta:=\Pi+\omega$, the \emph{Minkowski sum} of a polytope $\Pi$ and a cone $\omega$. The cone $\omega$ is called the \emph{tail} of $\Delta$. A polyhedron may have several Minkowski decompositions, but its tail is unique. A polyhedron whose tail is $\omega$ is called an \emph{$\omega$-tailed polyhedron}; if $\omega$ is pointed, $\Delta$ has finitely many vertices and is the Minkowski sum of their convex hull and $\omega$. The \emph{support function} of $\Delta$ is $h_{\Delta}:\omega^{\vee}\to \Q$, $h_{\Delta}(m):=\min_{v\in\Delta}\langle m,v\rangle$; it is positively homogeneous, concave $(h_{\Delta}(m)+h_{\Delta}(m')\leq h_{\Delta}(m+m'))$ and piecewise linear, and it determines $\Delta=\{v\in N_{\Q}\mid \langle m,v\rangle\geq h_{\Delta}(m)\textrm{ for all }m\in\omega^{\vee} \}$. Note that $h_{\omega}=0$.

For a lattice $N'$ with dual $M'$ and a pointed cone $\sigma\subset N'_{\Q}$, we denote by $X_{\sigma}:=\Spec(k[\sigma^{\vee}\cap M'])$ the associated affine toric variety, with torus $T_{N'}:=\Spec(k[M'])$. We freely use the orbit--cone correspondence \cite[Sections 1.3, 3.1 and 3.2]{CLS11}: the orbits of $X_{\sigma}$ correspond to the faces of $\sigma$, the (normal) toric prime divisors $D_{\varrho}$ to the rays $\varrho\in\sigma(1)$, with $X_{\sigma}\smallsetminus T_{N'}=\bigcup_{\varrho}D_{\varrho}$ and $\div(\chi^{m})=\sum_{\varrho}\langle m,n_{\varrho}\rangle D_{\varrho}$ for $m\in\sigma^{\vee}\cap M'$ ($n_{\varrho}$ the primitive generator of $\varrho$), and $(X_{\sigma})_{\chi^{m}}=X_{\varrho}$ if $m$ cuts out the face $\varrho$.

\subsection{Proper polyhedral divisors and normal $T$-varieties.} Let $Y$ be a normal variety over $k$ and $\omega\subset N_{\Q}$ be a pointed cone. A \emph{polyhedral divisor} over $Y$ with tail $\omega$ is a formal sum 
\[
\D:=\sum_{}\Delta_{D}\otimes D,
\]
where the sum runs over the prime Weil divisors of $Y$ and the coefficients $\Delta_{D}$ are $\omega$-tailed polyhedra, all but finitely many of which equal $\omega$. A polyhedral divisor is said to be \emph{proper} if $Y$ is semi-projective and it satisfies the following conditions:
    \begin{itemize}
        \item for all $m\in \omega^{\vee}\cap M$ the \emph{evaluation} 
            \[
            \D(m):=\sum h_{\Delta_{D}}(m)D=\sum\min_{v\in\Delta_D}\{\langle m,v\rangle\}D
            \]
        is a semiample $\Q$-Cartier divisor; and
        \item for all $m\in\mathrm{relint}(\omega^{\vee})\cap M$ the evaluation $\D(m)$ is a big $\Q$-Cartier divisor.
    \end{itemize}

Recall that a $\Q$-Cartier divisor $D\in\Div_{\Q}(Y)$ is \emph{semiample} if, for some positive integer $n$ such that $nD$ is Cartier, the open subsets $Y_{f}$, with $f\in \H^{0}(Y,\mathscr{O}(nD))$, cover $Y$. It is \emph{big} if, for some such $n$, there exists a section $f\in \H^{0}(Y,\mathscr{O}(nD))$ for which $Y_{f}$ is nonempty and affine.

Let $\D$ be a proper polyhedral divisor (hereafter referred to as a pp-divisor) over a normal semi-projective variety $Y$. By the concavity of the support functions, the evaluations $\D(m)$ satisfy a convexity property: for any pair $m$ and $m'$ in $\omega^{\vee}\cap M$, we have
\begin{equation}\label{equation: convexity}
    \D(m)+\D(m')\leq \D(m+m').
\end{equation}
This property allows us to construct the following $M$-graded algebras: 
\[\mathscr{A}(Y,\D):=\bigoplus_{m\in\omega^{\vee}\cap M}\mathscr{O}_{Y}(\D(m))\chi^{m} \quad\textrm{and}\quad A[Y,\D]:=\bigoplus_{m\in\omega^{\vee}\cap M}\H^{0}(Y,\mathscr{O}_{Y}(\D(m)))\chi^{m},\]
where the former is a quasi-coherent $\mathscr{O}_{Y}$-algebra and the latter a $k$-algebra. Here $\mathscr{O}_{Y}(\D(m))$ denotes $\mathscr{O}_{Y}(\lfloor\D(m)\rfloor)$, that is, for every open subset $V\subset Y$, 
\[\H^{0}(V,\mathscr{O}_{Y}(\D(m)))=\left\{f\in k(Y)^{*} \mid \left(\mathrm{div}(f)+\left\lfloor\D(m)\right\rfloor\right)|_{V}\geq 0\right\}\cup \{0\}.\]
Both algebras are contained in the group algebra $k(Y)[M]$, and the stalk of $\mathscr{A}(Y,\D)$ at the generic point of $Y$ is $k(Y)[\omega^{\vee}\cap M]$, with fraction field $k(Y)(M)$.

By \cite[Theorem 3.1]{AH06} and \cite[Proposition 4.17]{MN25a}, the scheme $X(\D):=\Spec(A[Y,\D])$ is an affine normal $T$-variety and $\tilde{X}(\D):=\Spec_{Y}(\mathscr{A}(Y,\D))$, the relative spectrum of the sheaf $\mathscr{A}(Y,\D)$, is a normal $T$-variety, the $T$-actions being given by the $M$-gradings. Moreover, by the same results, the affine normal $T$-variety $X(\D)$ is the affinisation of $\tilde{X}(\D)$, and the affinisation morphism $r_{\D}:\tilde{X}(\D)\to X(\D)$ is $T$-equivariant, proper and birational. These data fit into the following commutative diagram:
\[\xymatrix{ \tilde{X}(\D) \ar[r]^{r_{\D}} \ar[d]_{\pi} & X(\D) \ar[d]^{\pi_{0}} \\ Y \ar[r]_{r_{Y}} & \aff{Y}, }\]
where $\pi=\pi_{\D}$ is a $T$-invariant affine morphism. The function fields of $\tilde{X}(\D)$ and $X(\D)$ are both identified with $k(Y)(M)$, $\pi$ corresponding to the inclusion $k(Y)\subset k(Y)(M)$. In particular, $\dim(\tilde{X}(\D))=\dim(X(\D))=\dim(Y)+\rk(N)$.

\begin{remark}\label{remark: nonaffine base}
If $Y$ is not affine, then $\tilde{X}(\D)$ is not affine, and hence $\tilde{X}(\D)\not\cong X(\D)$. Indeed, $\pi$ is affine and $\mathscr{O}_{Y}=\mathscr{A}(Y,\D)_{0}$ is a direct summand of $\pi_{*}\mathscr{O}_{\tilde{X}(\D)}$. Hence, every quasi-coherent $\mathscr{O}_{Y}$-module $\mathscr{F}$ is a direct summand of $\pi_{*}\pi^{*}\mathscr{F}$, and $\H^{1}(Y,\pi_{*}\pi^{*}\mathscr{F})=\H^{1}(\tilde{X}(\D),\pi^{*}\mathscr{F})$ \cite[Exercise III.4.1]{Har77}. If $\tilde{X}(\D)$ were affine, the latter group would vanish \cite[Theorem III.3.5]{Har77}, so $\H^{1}(Y,\mathscr{F})=0$ for every quasi-coherent $\mathscr{O}_{Y}$-module $\mathscr{F}$, and $Y$ would be affine by Serre's criterion \cite[Theorem III.3.7]{Har77}. 
\end{remark}

Any affine normal $T$-variety arises from a pp-divisor by \cite[Theorem 3.4]{AH06} and \cite[Theorem 1.2 (ii)]{MN25a}. Comparing dimensions, the complexity of $X(\D)$ equals $\dim(Y)$.

A resolution of singularities of $\tilde{X}(\D)$ induces a resolution of $X(\D)$, since the affinisation morphism $r_{\D}:\tilde{X}(\D)\to X(\D)$ is proper and birational.

\subsection{Torus invariant divisors.} A $T$-invariant prime divisor $E$ of $\tilde{X}(\D)$ is called \emph{horizontal} if $\pi(E)$ is dense in $Y$ and \emph{vertical} otherwise. The general $T$-orbit of a horizontal (resp. vertical) prime divisor has dimension $\dim(T)-1$ (resp. $\dim(T)$). Over an algebraically closed field of characteristic zero, Petersen and S\"u\ss\ \cite[Proposition 3.13]{PS11} established bijections between the rays $\varrho\in\omega(1)$ and the horizontal prime divisors $\tilde{E}_{\varrho}$, and between the pairs $(D,v)$, where $D$ is a prime divisor of $Y$ and $v$ is a vertex of $\Delta_{D}$, and the vertical prime divisors $\tilde{E}_{(D,v)}$ of $\tilde{X}(\D)$. We do not use these bijections: over a smooth base with simple normal crossing support, the invariant prime divisors in the boundaries we consider are described directly by \cref{proposition: snc pp div and toroidal embeddings}, and the horizontal ones are defined intrinsically in \cref{Section: snc ppdivisors}.

\section{Resolution via toroidal embeddings for affine normal \texorpdfstring{$T$}{T}-varieties}\label{Section: equivariant resolution}

We first recall the definition of a \emph{toroidal embedding}, a notion developed in \cite{KKMSD73}. Then, we introduce the notion of \emph{simple normal crossing pp-divisor} and finally we prove \cref{theorem: main theorem affine} for affine varieties. In what follows, $X$ always stands for a normal variety.

\subsection{Toroidal embeddings.}\label{Section: toroidal embedding and equivariant resolution} Assume first that $X$ is a $T$-toric variety associated with a fan $\Sigma$ in $N_{\Q}$. By \cite[Theorem 11, p.~34]{KKMSD73}, there exists a $T$-invariant sheaf of ideals $\mathscr{I}\subset\mathscr{O}_{X}$ such that the blow-up along $\mathscr{I}$, denoted by $B_{\mathscr{I}}(X)$, is a $T$-equivariant resolution of singularities. It is therefore natural to ask whether normal varieties that \emph{locally} resemble toric varieties also admit desingularizations. Chapters I and II of \cite{KKMSD73}, which contain all the results used below, are written over an algebraically closed field of arbitrary characteristic.

\begin{definition}[\cite{KKMSD73}, Definition 1, p.~54]
    Let $X$ be a normal variety of dimension $n$ and $U\subset X$ be a smooth dense Zariski open subset. We say that $U\subset X$ is a \emph{toroidal embedding} if for every closed point $x\in X$, there exists an $n$-dimensional torus $T'$, an affine $T'$-toric variety $X_{\tau}$ associated with a pointed cone $\tau$, a closed point $t\in X_{\tau}$ and an isomorphism of complete local $k$-algebras
	\[\hat{\mathscr{O}}_{X,x}\cong \hat{\mathscr{O}}_{X_{\tau},t}\]
such that the ideal in $\hat{\mathscr{O}}_{X,x}$ generated by the ideal of $X\smallsetminus U$ (with its reduced structure) corresponds under this isomorphism to the ideal in $\hat{\mathscr{O}}_{X_{\tau},t}$ generated by the ideal of $X_{\tau}\smallsetminus T'$.
\end{definition}

Note that, in the previous definition, the torus $T'$ depends on the point $x\in X$. In particular, $X$ need not be a $T$-variety.

A Zariski open subvariety $U$ of a normal variety $X$ is said to be \emph{without self-intersection} if the irreducible components $E_{i}$ of $X\smallsetminus U$ are normal.

\begin{proposition}[\cite{KKMSD73}, Theorem 11*, p.~94]\label{proposition: Theorem 11* KKMSD73}
For any toroidal embedding $U\subset X$ without self-intersection, there exists a canonical coherent sheaf of ideals $\mathscr{I}\subset\mathscr{O}_{X}$ such that the normalization $\tilde{B}_{\mathscr{I}}(X)$ of the blow-up $B_{\mathscr{I}}(X)$ of $X$ along $\mathscr{I}$ is non-singular. In particular, $\Psi:\tilde{B}_{\mathscr{I}}(X)\to X$ is a projective resolution of singularities.
\end{proposition}

We have stated the conclusion for the normalized blow-up. If $B_{\mathscr{I}}(X)$ itself is non-singular, as one may also read \cite[Theorem 11*, p.~94]{KKMSD73}, then it coincides with its normalization, so the statement holds under either interpretation. The morphism $\Psi$ is projective, since blow-ups are projective and the normalization morphism is finite. It is also birational, since $\mathscr{I}\neq 0$.

We say that a toroidal embedding $U\subset X$ \emph{admits Zariski-local étale charts} if every closed point of $X$ has an open neighbourhood $V$ together with an étale morphism $\varphi: V\to X_{\tau}$ to an affine toric variety with dense torus $T'$ such that $\varphi^{-1}(X_{\tau}\smallsetminus T')=V\smallsetminus U$ as reduced closed subschemes. This is the notion of a \emph{strict toroidal embedding} of \cite[Definition 4.4.1]{Wlo20}. The notions of a toroidal embedding without self-intersection and a strict toroidal embedding are equivalent by \cite[Remark 4.4.4]{Wlo20}.

Suppose now that a normal $T$-variety $X$ admits a \emph{$T$-equivariant strict toroidal embedding} $U\subset X$, meaning that $U$ is $T$-stable. We use W{\l}odarczyk's canonical resolution in the following form. By \cite[Theorem 7.4.1]{Wlo20} there exists a canonical projective birational morphism $\Psi:X^{\mathrm{can}}\to X$, where $X^{\mathrm{can}}$ is smooth and $\Psi$ is an isomorphism over the non-singular locus of $X$, in particular over $U$. Moreover, the resolution is functorial with respect to smooth morphisms \cite[Definition 4.11.6]{Wlo20} and field extensions. Explicitly, for a smooth morphism $\varphi:X'\to X$, with the induced strict toroidal embedding $\varphi^{-1}(U)\subset X'$, the canonical resolution of $X'$ is the base change $X'\times_{X}X^{\mathrm{can}}\to X'$. Since $T$ is smooth and preserves $U$, both the action map $T\times X\to X$ and the projection $T\times X\to X$ are smooth and have inverse image $T\times U$. Functoriality therefore yields a $T$-action on $X^{\mathrm{can}}$ for which $\Psi$ is equivariant.

\begin{theorem}\label{theorem: equivariant resolution toroidal}
Let $X$ be a normal $T$-variety and $U\subset X$ a $T$-stable strict toroidal embedding. Then the $T$-action on $X$ lifts uniquely to a $T$-action on $X^{\mathrm{can}}$ for which $\Psi:X^{\mathrm{can}}\to X$ is a $T$-equivariant projective resolution of singularities, which is an isomorphism over $U$.
\end{theorem}

\begin{proof}
    Put $Z:=X\smallsetminus U$, a $T$-stable closed subset. Let $a,\mathrm{pr}:T\times X\to X$ be the action and the second projection, respectively. Both morphisms are smooth since $\mathrm{pr}$ is a base change of $T\to\Spec(k)$, and $a=\mathrm{pr}\circ\theta$, where $\theta:(t,x)\mapsto(t,t\cdot x)$ is an automorphism of $T\times X$. Moreover, since $Z$ is $T$-stable, both satisfy $a^{-1}(Z)=T\times Z=\mathrm{pr}^{-1}(Z)$. Moreover, $T\times U$ is also a strict toroidal embedding. Indeed, the irreducible components of $T\times Z$ are products $T\times E$, for $E$ running through the irreducible components of $Z$. These products are normal, since they are smooth over the normal varieties $E$. If $\varphi:V\to X_{\tau}$ is a chart of $X$, then $\mathrm{id}\times\varphi:T\times V\to T\times X_{\tau}$ is an étale morphism to the affine toric variety associated with the cone $\{0\}\times\tau$ for the torus $T\times T'$. The inverse image of its boundary $T\times (X_{\tau}\smallsetminus T')$ is $T\times(V\smallsetminus U)$.

    By functoriality of the canonical resolution with respect to the smooth morphisms $\mathrm{pr}$ and $a$, the canonical resolution of $T\times X$ is identified with each of
    \[
    (T\times X)\times_{\mathrm{pr},X}X^{\mathrm{can}}\cong T\times X^{\mathrm{can}}\quad\mathrm{ and }\quad(T\times X)\times_{a,X}X^{\mathrm{can}}
    \]
    over $T\times X$. The resulting isomorphism restricts to the identity over $T\times U$, where both sides are identified with $T\times U$ via their structure morphisms. Composing this isomorphism with the projection to $X^{\mathrm{can}}$ gives a morphism 
    \[
    \alpha:T\times X^{\mathrm{can}}\to X^{\mathrm{can}}\quad\textrm{such that}\quad \Psi\circ\alpha=a\circ(\mathrm{id}_{T}\times \Psi).
    \]
    Since $U$ is $T$-stable and $\Psi$ restricts to an isomorphism $\Psi^{-1}(U)\to U$, the restriction of $\alpha$ to $T\times \Psi^{-1}(U)$ is the $T$-action on $U$ transported via this isomorphism. Hence the identities
    \[
    \alpha(1,w)=w \quad\textrm{and}\quad\alpha(t,\alpha(t',w))=\alpha(tt',w)
    \]
    hold on the dense open subsets $\Psi^{-1}(U)$ and $T\times T\times \Psi^{-1}(U)$, respectively. Since the schemes $X^{\mathrm{can}}$ and $T\times T\times X^{\mathrm{can}}$ are reduced and $X^{\mathrm{can}}$ is separated, these identities hold everywhere \cite[Exercise II.4.2]{Har77}. Thus, $\alpha$ defines an action of $T$ on $X^{\mathrm{can}}$ making $\Psi$ equivariant. The action is effective because it is effective on $\Psi^{-1}(U)\cong U$, and it is unique by the same density argument.
\end{proof}

\subsection{Simple normal crossing pp-divisors.}\label{Section: snc ppdivisors}


Let $Y$ be a smooth variety, not necessarily semi-projective, and $\D=\sum\Delta_D\otimes D$ be a polyhedral divisor on $Y$ with pointed tail $\omega$, not necessarily proper. Since $Y$ is smooth, all $\D(m)$ are $\Q$-Cartier and $\mathscr{A}(Y,\D)$ is defined as before. We define the \emph{support} of $\D$ as the reduced effective divisor
\[
\Supp(\D):=\sum_{\Delta_D\neq\omega}D.
\]
Let $B$ be a reduced effective divisor on $Y$ with $\Supp(\D)\leq B$ and put $Y^{\circ}:=Y\smallsetminus B$. All coefficients of $\D$ on $Y^{\circ}$ equal $\omega$, so $\D(m)|_{Y^{\circ}}=0$ and $\mathscr{A}(Y,\D)|_{Y^{\circ}}=\mathscr{O}_{Y^{\circ}}\otimes_k k[\omega^{\vee}\cap M]$. Hence, whenever $\mathscr{A}(Y,\D)$ is of finite type and $\tilde{X}(\D):=\Spec_Y(\mathscr{A}(Y,\D))$ with projection $\pi:\tilde{X}(\D)\to Y$, there is a canonical $T$-equivariant isomorphism
\begin{equation}\label{equation: trivial boundary complement}
\pi^{-1}(Y^{\circ})\cong Y^{\circ}\times X_{\omega}.
\end{equation}
For $\varrho\in\omega(1)$, let $\tilde{E}_{\varrho}$ be the closure in $\tilde{X}(\D)$ of $Y^{\circ}\times D_{\varrho}$. These are pairwise distinct horizontal $T$-invariant prime divisors. We set
\[
Z_B(\D):=\pi^{-1}(B)\cup\left(\bigcup_{\varrho\in\omega(1)}\tilde{E}_{\varrho}\right)
\quad\text{and}\quad
U_B(\D):=\tilde{X}(\D)\smallsetminus Z_B(\D)=Y^{\circ}\times T,
\]
where the last equality follows from \eqref{equation: trivial boundary complement}, since $X_{\omega}\smallsetminus T=\bigcup_{\varrho}D_{\varrho}$. Thus, $Z_B(\D)$ is a $T$-stable closed subset. Note that a closed point $x\in\tilde{X}$ lies in $U_B(\D)$ if and only if $\pi(x)\notin B$ and $\dim(T\cdot x)=\dim(T)$, since if $\pi(x)=y\notin B$, then $x=(y,z)$ with $z\in X_{\omega}$ under \eqref{equation: trivial boundary complement}, $T\cdot x=\{y\}\times T\cdot z$, and $\dim(T\cdot z)=\dim(T)$ if and only if $z\in T$. When $B=\Supp(\D)$, we write $Z(\D):=Z_B(\D)$ and $U(\D):=U_B(\D)$.

The following result states that $U_{B}(\D)\subset\tilde{X}(\D)$ is a strict toroidal embedding (or without self-intersection) whenever $\D$ is a polyhedral divisor over a smooth variety $Y$ and $B$ is a simple normal crossing divisor. Its proof gives a characteristic-free local description of $\tilde{X}(\D)$, in which the cone $\sigma$ is the Cayley cone of \cite[Example 2.5]{LS13}. The freedom in the choice of $B$ will be used in \cref{Section: resolution general case}.

\begin{proposition}\label{proposition: snc pp div and toroidal embeddings}
Let $Y$ be a smooth variety of dimension $n$, $\D$ a polyhedral divisor on $Y$ with pointed tail $\omega$, and $B$ a simple normal crossing divisor on $Y$ with $\Supp(\D)\leq B$.
\begin{enumerate}
\item \label{proposition: snc pp div and toroidal embeddings part i} $\mathscr{A}(Y,\D)$ is a finitely generated $\mathscr{O}_Y$-algebra and $\tilde{X}(\D)=\Spec_Y(\mathscr{A}(Y,\D))$ is a normal variety of dimension $n+\rk(N)$ with an effective $T$-action and a $T$-invariant affine morphism $\pi:\tilde{X}(\D)\to Y$.
\item \label{proposition: snc pp div and toroidal embeddings part ii} Let $y\in Y$ be a closed point, $D_1,\dots,D_r$ the components of $B$ through $y$, $\Delta_i:=\Delta_{D_i}$, and let $\sigma\subset N_{\Q}\times\Q^r$ be the cone generated by $\omega\times\{0\}$ and $\Delta_i\times\{e_i\}$, $1\leq i\leq r$. Then $\sigma$ is pointed, and there exist an affine open neighbourhood $V$ of $y$ and an étale morphism $q:\pi^{-1}(V)\to X_{\sigma}\times\Gm^{n-r}$ with
\[
Z_B(\D)\cap\pi^{-1}(V)=q^{-1}((X_{\sigma}\smallsetminus T_{N\times\Z^r})\times\Gm^{n-r})
\]
as reduced closed subschemes.
\item \label{proposition: snc pp div and toroidal embeddings part iii} $U_B(\D)\subset\tilde{X}(\D)$ is a toroidal embedding without self intersection.
\end{enumerate}
\end{proposition}

\begin{proof}
\emph{Step 1: coordinates and the cone $\sigma$.} Choose a regular system of parameters $y_1,\dots,y_n$ of $\mathscr{O}_{Y,y}$ with $y_i$ a local equation of $D_i$ for $i\leq r$, and put $u_j:=1+y_j$ for $j>r$. After shrinking, we obtain an affine open neighbourhood $V$ of $y$ such that $y_i,u_j\in\mathscr{O}(V)$, the $u_j$ are units, $\div(y_i)|_V=D_i|_V$, $V$ meets no other component of $B$, and $\varphi:=(y_1,\dots,y_r,u_{r+1},\dots,u_n):V\to\A^r\times\Gm^{n-r}$ is étale. Indeed, $dy_1,\dots,dy_n$ form a basis of $\Omega_Y\otimes k(y)$, so $\Omega_{V/\A^n}$ vanishes near $y$, and a morphism of smooth varieties of the same dimension with vanishing relative differentials is étale \cite[Proposition III.10.4]{Har77}.

Let $\sigma$ be the cone defined above, which depends on $\D$. Every element of $\sigma$ has the form $(w+\sum_i\lambda_iv_i,\lambda)$ with $w\in\omega$, $\lambda\in\Q^r_{\geq0}$ and $v_i\in\Delta_i$, by convexity of the $\Delta_i$. Let us prove that $\sigma$ is pointed. If $\pm x\in\sigma$, then $\lambda=0$ and $\pm x\in\omega\times\{0\}$, so $x=0$. Writing $h_i:=h_{\Delta_i}$, a vector $(m,a)\in M_{\Q}\times\Q^r$ lies in $\sigma^{\vee}$ if and only if $m\in\omega^{\vee}$ and $a_i\geq-h_i(m)$ for all $i$. The face of $\sigma$ cut out by $(0,\sum_ie_i)$ is $\omega\times\{0\}$, and the face cut out by $(0,\sum_{j\neq i}e_j)$ is the cone over $(\Delta_i\times\{1\})\cup(\omega\times\{0\})$. Since every ray of $\sigma$ is spanned by a generator, the rays of $\sigma$ are the $\varrho\times\{0\}$, $\varrho\in\omega(1)$, and the rays spanned by $(v,e_i)$, $v$ a vertex of $\Delta_i$.

\emph{Step 2: the local model.} For $m\in\omega^{\vee}\cap M$, we have $\lfloor\D(m)\rfloor|_V=\div(\prod_i y_i^{\lfloor h_i(m)\rfloor})|_V$. Hence, since $V$ is normal, $\H^0(V,\mathscr{O}_Y(\D(m)))=\mathscr{O}(V)\prod_i y_i^{-\lfloor h_i(m)\rfloor}$. On the other hand, with $z_i:=\chi^{(0,e_i)}$ and $k[z]:=k[z_1,\dots,z_r]$, the description of $\sigma^{\vee}$ gives
\[
k[\sigma^{\vee}\cap(M\times\Z^r)]=\bigoplus_{m\in\omega^{\vee}\cap M}k[z]\prod_i z_i^{-\lfloor h_i(m)\rfloor}\chi^m,
\]
each summand being free of rank one over $k[z]$. Hence the homomorphism of $M$-graded $\mathscr{O}(V)$-algebras
\[
\mathscr{O}(V)\otimes_{k[z]}k[\sigma^{\vee}\cap(M\times\Z^r)]\longrightarrow\H^0(V,\mathscr{A}(Y,\D)),
\qquad f\otimes z^a\chi^m\longmapsto fy^a\chi^m,
\]
where $z_i\mapsto y_i$, is an isomorphism. In each degree, it maps a free $\mathscr{O}(V)$-module of rank one onto another, sending a generator to a generator. Thus, $\mathscr{A}(Y,\D)|_V$ is of finite type and
\begin{equation}\label{equation: local model}
\pi^{-1}(V)\cong V\times_{\A^r}X_{\sigma}
=V\times_{\A^r\times\Gm^{n-r}}(X_{\sigma}\times\Gm^{n-r}),
\end{equation}
where $V\to\A^r$ is $(y_i)$ and $X_{\sigma}\to\A^r$ is $(z_i)$. Let $q:\pi^{-1}(V)\to X_{\sigma}\times\Gm^{n-r}$ be the base change of $\varphi$, which is étale, and let $p:\pi^{-1}(V)\to X_{\sigma}$ be its composition with the projection, which is smooth, hence flat and open.

\emph{Proof of \ref{proposition: snc pp div and toroidal embeddings part i}.} Since $y$ is arbitrary, $\mathscr{A}(Y,\D)$ is of finite type, and $\tilde{X}(\D)$ is of finite type, separated and integral ($\mathscr{A}(Y,\D)\subset k(Y)[M]$ is a domain). By \eqref{equation: local model}, $\pi^{-1}(V)$ is smooth over the normal variety $X_{\sigma}$, hence normal by \cite[Theorem 23.9 and its Corollary]{Mat86}, of dimension $\rk N+r+(n-r)$. The action induced by the grading is effective because $\tilde{X}(\D)\supset Y^{\circ}\times T$.

\emph{Proof of \ref{proposition: snc pp div and toroidal embeddings part ii}.} Put $\partial:=X_{\sigma}\smallsetminus T_{N\times\Z^r}=\bigcup_{\varsigma\in\sigma(1)}D_{\varsigma}$. First,
\[
\pi^{-1}(B\cap V)=\pi^{-1}(V(y_1\cdots y_r))=p^{-1}(V(z_1\cdots z_r)),
\]
and $\div(z_i)=\sum_{\varsigma}\langle(0,e_i),n_{\varsigma}\rangle D_{\varsigma}$ has a positive coefficient exactly along the $D_{\varsigma}$ with $\varsigma$ spanned by some $(v,e_i)$. Then, $V(z_1\cdots z_r)$ is the union of these $D_{\varsigma}$. Second, $X_{\sigma}\smallsetminus V(z_1\cdots z_r)$ is the toric variety $X_{\omega}\times\Gm^r$ of the face $\omega\times\{0\}$, and over $V\smallsetminus B$ the isomorphism \eqref{equation: local model} restricts to \eqref{equation: trivial boundary complement}. Hence $(V\smallsetminus B)\times(X_{\omega}\smallsetminus T)=p^{-1}(S_0)$ with $S_0:=\bigcup_{\varrho}D_{\varrho\times\{0\}}\smallsetminus V(z_1\cdots z_r)$. Each $D_{\varrho\times\{0\}}$ is irreducible and meets $X_{\sigma}\smallsetminus V(z_1\cdots z_r)$ (the $z_i$ do not vanish on the orbit of $\varrho\times\{0\}$, since $(0,e_i)\in(\varrho\times\{0\})^{\perp}$), so $\overline{S_0}=\bigcup_{\varrho}D_{\varrho\times\{0\}}$, and, since $p$ is open, $\overline{p^{-1}(S_0)}=p^{-1}(\overline{S_0})$. Since $\tilde{E}_{\varrho}\cap\pi^{-1}(V)$ is the closure of $(V\smallsetminus B)\times D_{\varrho}$, we obtain
\[
Z_B(\D)\cap\pi^{-1}(V)=p^{-1}(\partial)=q^{-1}(\partial\times\Gm^{n-r})
\]
set-theoretically. As $q$ is étale, $q^{-1}(\partial\times\Gm^{n-r})$ is reduced \cite[Corollary to Theorem 23.9]{Mat86}, and the equality holds as reduced closed subschemes.

\emph{Proof of \ref{proposition: snc pp div and toroidal embeddings part iii}.} Let $x\in\tilde{X}(\D)$ be a closed point. By part \ref{proposition: snc pp div and toroidal embeddings part ii}, there exists an open neighbourhood $V\subset Y$ of $\pi(x)$ and an \'etale morphism $q:\pi^{-1}(V)\to X_{\tau}$ such that \[
Z_B(\D)\cap\pi^{-1}(V)=q^{-1}(X_{\tau}\smallsetminus T')
\]
as reduced closed subschemes, where $X_{\tau}$ is a $T'$-toric variety. Since $x$ is arbitrary, these morphisms provide Zariski-local \'etale charts for $U_{B}(\D)\subset \tilde{X}(\D)$. Thus, it is a strict toroidal embedding, or equivalently a toroidal embedding without self-intersection.

\end{proof}

\begin{definition}\label{definition: simple normal crossing}
A pp-divisor $\D$ over a normal semi-projective variety $Y$ is called \emph{simple normal crossing (SNC)} if $Y$ is smooth and $\Supp(\D)$ is a simple normal crossing divisor.
\end{definition}

\begin{corollary}\label{corollary: snc pp div and toroidal embeddings}
Let $\D$ be a pp-divisor over a smooth semi-projective variety $Y$ and $B$ be a simple normal crossing divisor on $Y$ such that $\mathrm{Supp}(\D)\leq B$, then $U_{B}(\D)\subset\tilde{X}(\D)$ is a toroidal embedding without self-intersection. In particular, if $\D$ is a simple normal crossing pp-divisor, then $U(\D)\subset\tilde{X}(\D)$ is a toroidal embedding without self-intersection.
\end{corollary}

\begin{proof}
    This is direct consequence of \cref{proposition: snc pp div and toroidal embeddings}\ref{proposition: snc pp div and toroidal embeddings part iii}.
\end{proof}

\begin{remark}\label{remark: normal crossing support}
If $\Supp(\D)$ only has normal crossings (i.e. the coordinates of Step 1 exist étale locally), the same computation over an étale neighbourhood shows that $U(\D)\subset\tilde{X}(\D)$ is a toroidal embedding. The simple normal crossing hypothesis gives Zariski coordinates defining the global components of $B$, which ensures that the embedding is without self-intersection (compare \cite[Proposition 2.6 and Remark 2.8]{LS13} in characteristic zero). The prime divisors with $\Delta_D=\omega$ play no role in the toroidal structure, since over them $\tilde{X}(\D)$ is locally $Y\times X_{\omega}$. For this reason, unlike \cite[Definition 2.2]{LS13}, we do not include in $\Supp(\D)$ the divisors with trivial coefficient and no big vertex and any enlargement of the boundary is covered by the choice of $B$.
\end{remark}

\subsection{Pull-backs of pp-divisors and log resolutions.}\label{Section: pullbacks}
Let $X$ be an affine normal $T$-variety. By \cite[Theorem 3.4]{AH06} and \cite[Theorem 1.2(ii)]{MN25a}, there exists a pp-divisor $\D$ over a normal semi-projective variety $Y$ with tail cone $\omega\subset N_{\Q}$ such that $X\cong X(\D)$ as $T$-varieties. This algebro-geometric combinatorial description is far from unique. In fact, two non-isomorphic pp-divisors $\D$ and $\D'$ can give rise to isomorphic affine normal $T$-varieties, as shown in \cref{lemma: birational modification of a pp-div} below.

Let $\psi:Y'\to Y$ be a proper birational morphism of varieties with $Y$ normal, and let $W\subset Y$ be the complement of the largest open subset over which $\psi$ is an isomorphism. Then $W$ has codimension at least two: if $\eta$ is the generic point of a prime divisor of $Y$, then $Y'\times_Y\Spec(\mathscr{O}_{Y,\eta})$ is integral, proper and birational over the spectrum of a discrete valuation ring, hence one-dimensional with finite closed fibre. By upper semicontinuity of fibre dimension \cite[Exercise II.3.22]{Har77}, $\psi$ is quasi-finite, hence finite \cite[Exercise III.11.2]{Har77}, over a neighbourhood of $\eta$, and a finite birational morphism onto a normal variety is an isomorphism. We call $\operatorname{Exc}(\psi):=\psi^{-1}(W)$ the \emph{exceptional set} of $\psi$. A prime divisor $E\not\subset\operatorname{Exc}(\psi)$ of $Y'$ meets $\psi^{-1}(Y\smallsetminus W)\cong Y\smallsetminus W$, hence is the strict transform of the prime divisor $\psi(E)$, and for every $\Q$-Cartier divisor $F$ on $Y$ the coefficient of $E$ in $\psi^{*}F$ equals that of $\psi(E)$ in $F$.

\begin{lemma}\label{lemma: birational modification of a pp-div}
Let $\D$ be a pp-divisor over a normal semi-projective variety $Y$ with tail $\omega$, and let $\psi:Y'\to Y$ be a proper birational morphism with $Y'$ normal.
\begin{enumerate}
\item \label{lemma: birational modification of a pp-div part i} There is a unique polyhedral divisor $\psi^*\D=\sum_E\Delta'_{E}\otimes F$ on $Y'$ with tail $\omega$ whose evaluations are the pullbacks $(\psi^*\D)(m)=\psi^*(\D(m))$ for all $m\in\omega^{\vee}\cap M$.
\item \label{lemma: birational modification of a pp-div part ii}If $E\not\subset\mathrm{Exc}(\psi)$, then $\Delta'_{E}=\Delta_{\psi(E)}$. In particular, every component of $\Supp(\psi^*\D)$ is exceptional or lies in $\psi^{-1}(\Supp(\D))$.
\item \label{lemma: birational modification of a pp-div part iii} In $k(Y')=k(Y)$, we have
\[
\H^0(Y',\mathscr{O}_{Y'}(\psi^*\D(m)))=\H^0(Y,\mathscr{O}_Y(\D(m)));
\]
hence $A[Y',\psi^*\D]=A[Y,\D]$.
\item \label{lemma: birational modification of a pp-div part iv} If $\psi$ is projective, then $Y'$ is semi-projective, $\psi^*\D$ is a pp-divisor, and $X(\psi^*\D)=X(\D)$ as $T$-varieties.
\end{enumerate}
\end{lemma}

\begin{proof}
\emph{Proof of \ref{lemma: birational modification of a pp-div part i}} Let $E$ be a prime divisor of $Y'$ and $m\in\omega^{\vee}\cap M$. Let $h_{E}(m)$ be the coefficient of $E$ in $\psi^*(\D(m))$, where $\D(m)=\sum h_{\Delta_{D}}(m)D$. The function $h_{E}$ is positively homogeneous and concave, since $\D(m+m')-\D(m)-\D(m')$ is effective and $\Q$-Cartier by \eqref{equation: convexity}, and pull-back preserves effectivity. Let $\Lambda$ be the normal fan of the $\omega$-tailed polyhedron $P:=\sum_{D}\Delta_{D}$ (Minkowski sum over the finitely many $D$ with $\Delta_{D}\neq\omega$ or $P:=\omega$ if there are none). Its maximal cones $\lambda$ are full-dimensional, cover $\omega^{\vee}$, and every $h_{\Delta_{D}}$ is linear on each of them. Thus $m\mapsto\D(m)$ is linear on $\lambda$, and so is $h_{E}$, say $h_{E}|_{\lambda}=\langle\cdot,v_{\lambda}\rangle$ with $v_{\lambda}\in N_{\Q}$. For $m\in\omega^{\vee}$ and $m_{0}$ interior to $\lambda$, concavity gives
\[
h_{E}(tm+(1-t)m_0)\geq t\cdot h_{E}(m)+(1-t)\cdot h_{E}(m_0),
\]
and for $0<t\ll1$ the left-hand side is $\langle tm+(1-t)m_0,v_{\lambda}\rangle$. Hence $h_{E}(m)\leq\langle m,v_{\lambda}\rangle$, with equality on $\lambda$. Thus $h_{E}=\min_{\lambda}\langle\cdot,v_{\lambda}\rangle$ is the support function of $\Delta'_{E}:=\operatorname{conv}\{v_{\lambda}\}+\omega$, which is determined by $h_{E}$. Besides, $\Delta'_{E}=\omega$ for all but finitely many $E$ by \ref{lemma: birational modification of a pp-div part ii}, since there are finitely many exceptional prime divisors.

\emph{Proof of \ref{lemma: birational modification of a pp-div part ii}} By the discussion above the lemma, we have that $E$ is the strict transform of the prime divior $\psi(E)$ and hence $h_{E}=h_{\Delta_{\psi(E)}}$.

\emph{Proof of \ref{lemma: birational modification of a pp-div part iii}}  Set $F:=\D(m)$ and let $f\in k(Y)^{*}$. If $\div(f)+\lfloor F\rfloor\geq0$, then $\div(f)+F\geq0$ is $\Q$-Cartier, so $\div(f)+\psi^{*}F=\psi^{*}(\div(f)+F)\geq0$ and, since $\div(f)$ is integral, $\div(f)+\lfloor\psi^{*}F\rfloor\geq0$. Conversely, let $\psi_*$ be the push-forward of Weil divisors (killing exceptional primes). It preserves effectivity, $\psi_{*}\div_{Y'}(f)=\div_{Y}(f)$ and $\psi_{*}\lfloor\psi^{*}F\rfloor=\lfloor F\rfloor$ by \ref{lemma: birational modification of a pp-div part ii}, so $\div(f)+\lfloor\psi^{*}F\rfloor\geq0$ implies $\div(f)+\lfloor F\rfloor\geq0$.

\emph{Proof of \ref{lemma: birational modification of a pp-div part iv}} The variety $Y'$ is semi-projective by \cite[Lemma 4.11]{MN25a}. Pull-backs of semiample divisors are semiample, since $Y'_{\psi^*f}=\psi^{-1}(Y_f)$. For $m\in\operatorname{relint}(\omega^{\vee})\cap M$, choose $n\in \N$ such that $L:=n\D(m)$ is Cartier and $f\in\H^0(Y,\mathscr{O}_Y(L))$ with $Y_f$ non-empty and affine. Choose $h\in\mathscr{O}(Y_f)$ with $\emptyset\neq(Y_f)_h\subset Y\smallsetminus W$. By \cite[Lemma II.5.14]{Har77}, $f^sh$ extends to some $a\in\H^0(Y,\mathscr{O}_Y(sL))$, and $Y_{fa}=(Y_f)_h$. Then $Y'_{\psi^*(fa)}=\psi^{-1}(Y_{fa})\cong Y_{fa}$ is non-empty and affine, so $\psi^*\D(m)$ is big. The last assertion follows from \ref{lemma: birational modification of a pp-div part i}.
\end{proof}

In characteristic zero, \cref{lemma: birational modification of a pp-div}\ref{lemma: birational modification of a pp-div part iv} is \cite[Corollary 8.12]{AH06} (see also \cite[Lemma 2.1]{LS13}). Over arbitrary fields, a closely related statement is contained in \cite[Theorem 5.11]{MN25a}.

Let $\D$ be a pp-divisor over a normal semi-projective variety $Y$ with tail cone $\omega\subset N_{\Q}$ such that $X\cong X(\D)$ as $T$-varieties. By \cref{lemma: birational modification of a pp-div}, for any projective birational morphism $\psi:Y'\to Y$, the varieties $X(\D')$ and $X(\D)$ are isomorphic as $T$-varieties, where $\D':=\psi^*\D$. However, the $T$-varieties $\tilde{X}(\D)$ and $\tilde{X}(\D')$ need not be isomorphic. For instance, if $\D$ is the trivial pp-divisor on $Y=\A^2$ and $\psi$ is the blow-up of the origin, then $\tilde{X}(\D)=\A^2\times X_{\omega}$ is affine while $\tilde{X}(\D')=Y'\times X_{\omega}$ is not, by \cref{remark: nonaffine base}. The normal $T$-variety $\tilde{X}(\D)$ associated with one pp-divisor may admit no toroidal embedding as above, whereas the variety associated with a non-isomorphic pp-divisor may admit one, for instance when the latter pp-divisor is SNC. We cannot ensure the existence of an SNC pp-divisor for affine normal $T$-varieties in arbitrary complexity because of the smoothness condition on $Y$. In low complexity, however, this can be ensured by the following log resolutions.

\begin{proposition}\label{proposition: log resolution base}
Let $Y$ be a normal variety of dimension at most three and $B\subset Y$ a closed subset of codimension at least one, possibly empty. There exists a proper birational morphism $\psi:Y'\to Y$ from a smooth variety $Y'$ which is an isomorphism over $Y\smallsetminus(B\cup\Sing(Y))$ and such that $\psi^{-1}(B\cup\Sing(Y))$, with its reduced structure, is the support of a simple normal crossing divisor $B'$ on $Y'$ (possibly $B'=0$). In particular, $\operatorname{Exc}(\psi)\cup\psi^{-1}(B)\subset B'$. If $\dim(Y)\leq2$, then $\psi$ can be chosen projective.
\end{proposition}

\begin{proof}
If $\dim(Y)\leq1$, then $Y$ is a point or a normal curve, hence smooth, and $B$ is a finite set of points, which is a simple normal crossing divisor, and in such a case take $\psi=\mathrm{id}$. Assume $\dim(Y)\in\{2,3\}$. The variety $Y$ is a reduced, separated, quasi-excellent Noetherian scheme of dimension at most three, so the theorem of Cossart and Piltant \cite[Theorem 1.1]{CP19} gives a proper birational morphism $\pi_1:Y_1\to Y$ with $Y_1$ regular which is an isomorphism over $Y\smallsetminus\Sing(Y)$. The closed subset $X_1:=\pi_1^{-1}(B\cup\Sing(Y))$, with its reduced structure, has dimension at most two. By \cite[Corollary 0.4]{CJS09}, applied to the regular excellent scheme $Y_1$ and to $X_1$, there is a projective surjective morphism $\pi_2:Y'\to Y_1$ which is an isomorphism over $Y_1\smallsetminus X_1$ and such that $\pi_2^{-1}(X_1)_{\mathrm{red}}$ is a simple normal crossing divisor on $Y'$. Moreover, $Y'$ is regular, since $\pi_2$ is constructed as a composition of blow-ups of regular schemes along regular centres. The scheme $Y'$ is integral, since it is reduced, $\pi_2$ is an isomorphism over the dense open subset $Y_1\smallsetminus X_1$, and the divisor $\pi_2^{-1}(X_1)$ contains no irreducible component of $Y'$. Thus $\psi:=\pi_1\circ\pi_2$ is proper and birational. Besides, $Y'$ is smooth, $\psi$ is an isomorphism over $Y\smallsetminus(B\cup\Sing(Y))$ and $\psi^{-1}(B\cup\Sing(Y))=\pi_2^{-1}(X_1)$. Since the set $W$ of \cref{Section: pullbacks} is contained in $B\cup\Sing(Y)$, we obtain $\operatorname{Exc}(\psi)\subset B'$. For surfaces, one may instead use Lipman's resolution \cite{Lip78}, obtained by alternately normalizing and blowing up singular points, hence projective, followed by the embedded resolution of the curve $\pi_1^{-1}(B\cup\Sing(Y))$ in the smooth surface $Y_1$ by point blow-ups \cite[Theorem V.3.9]{Har77}. This yields a projective $\psi$.
\end{proof}

For threefolds in characteristic $p>5$, resolution and embedded resolution of surfaces go back to Abhyankar \cite{Abh66} (see also \cite{Cut09}). In dimension three, the morphism $\psi$ of \cref{proposition: log resolution base} is only asserted to be proper, which is why the following statement is conditional.

\begin{proposition}\label{proposition: existence of snc pp-divisors}
Let $X$ be an affine normal $T$-variety of complexity at most three, and let $\D$ be a pp-divisor over a normal semi-projective variety $Y$ such that $X\cong X(\D)$. If $\dim(Y)\leq2$, or if the morphism $\psi$ of \cref{proposition: log resolution base} for $(Y,\Supp(\D))$ can be chosen projective, then there exists an SNC pp-divisor $\D'$ such that $X\cong X(\D')$.
\end{proposition}

\begin{proof}
We have $\dim(Y)\leq3$, since $\dim(Y)$ is the complexity of $X$. If $\dim(Y)\leq1$, then $Y$ is already smooth and $\Supp(\D)$ is a finite set of points, so $\D$ is automatically an SNC pp-divisor. Otherwise, let $\psi:Y'\to Y$ and $B'$ be as in \cref{proposition: log resolution base} for $B:=\Supp(\D)$, with $\psi$ projective, as allowed by the assumption and, for surfaces, by the last assertion of that proposition. By \cref{lemma: birational modification of a pp-div}\ref{lemma: birational modification of a pp-div part ii}, the components of $\Supp(\psi^*\D)$ are among those of $B'$, so $\Supp(\psi^*\D)$ is a simple normal crossing divisor. By \cref{lemma: birational modification of a pp-div}\ref{lemma: birational modification of a pp-div part iv}, $\D':=\psi^*\D$ is an SNC pp-divisor over the smooth semi-projective variety $Y$ such that $X\cong X(\D')\cong X(\D)$ as $T$-varieties.
\end{proof}

The following proposition shows that, over a smooth base with simple normal crossing support, the pull-back along a merely proper birational morphism still yields a proper birational modification of $X(\D)$ with a toroidal embedding. Thus, \cref{theorem: main theorem affine} holds without a projectivity assumption.

\begin{proposition}\label{proposition: proper pullback modification}
Let $\D$ be a pp-divisor over a normal semi-projective variety $Y$ with tail $\omega$. Let $\psi:Y'\to Y$ be a proper birational morphism from a smooth variety $Y'$, and $B'$ a simple normal crossing divisor on $Y'$ with $\Supp(\psi^*\D)\leq B'$. Let $\tilde{X}':=\tilde{X}(\psi^*\D)$ with projection $\pi':\tilde{X}'\to Y'$. Then $\tilde{X}'$ is a normal $T$-variety, $\H^0(\tilde{X}',\mathscr{O}_{\tilde{X}'})=A[Y,\D]$, and the induced morphism $r':\tilde{X}'\to\Spec(\H^0(\tilde{X}',\mathscr{O}_{\tilde{X}'}))=X(\D)$ is $T$-equivariant, proper and birational. Moreover, $U_{B'}(\psi^*\D)\subset\tilde{X}'$ is a $T$-stable toroidal embedding without self-intersection.
\end{proposition}

\begin{proof}
By \cref{proposition: snc pp div and toroidal embeddings}, applied to $\psi^*\D$ and $B'$, the algebra $\mathscr{A}':=\mathscr{A}(Y',\psi^*\D)$ is of finite type, $\tilde{X}'$ is a normal $T$-variety and $U_{B'}(\psi^*\D)\subset\tilde{X}'$ is as claimed. By \cref{lemma: birational modification of a pp-div}\ref{lemma: birational modification of a pp-div part iii}, $\H^0(\tilde{X}',\mathscr{O}_{\tilde{X}'})=\H^0(Y',\mathscr{A}')=A[Y,\D]$, and $r'$ is $T$-equivariant since it respects the gradings. Put $\mathscr{A}:=\mathscr{A}(Y,\D)$. The argument of \cref{lemma: birational modification of a pp-div}\ref{lemma: birational modification of a pp-div part iii}, applied over open subsets $V\subset Y$, gives
\[
\H^0(V,\mathscr{O}_Y(\D(m)))\subset\H^0(\psi^{-1}(V),\mathscr{O}_{Y'}(\psi^*\D(m))),
\]
hence a homomorphism of graded $\mathscr{O}_{Y'}$-algebras $\psi^*\mathscr{A}\to\mathscr{A}'$ whose image we denote by $\mathscr{B}$. We claim that $\mathscr{A}'$ is integral over $\mathscr{B}$. Let $y'\in Y'$, let $V\subset Y$ be an affine open neighbourhood of $\psi(y')$ on which $l\D(m)$ is Cartier and equal to $\div(g)$, and $V''\subset\psi^{-1}(V)$ an affine open neighbourhood of $y'$. Then $g^{-1}\chi^{lm}\in\H^0(V,\mathscr{A})$. For $a=f\chi^m\in\H^0(V'',\mathscr{A}')$, we have $\div(f)+\psi^*\D(m)\geq0$ on $V''$, so $\div(f^lg)\geq0$ and $f^lg\in\mathscr{O}(V'')$ by normality. Thus $a^l=(f^lg)(g^{-1}\chi^{lm})\in\H^0(V'',\mathscr{B})$. Since $\mathscr{A}'$ is of finite type, it is a finite $\mathscr{B}$-module, so the natural morphism
\[
\begin{split}
\Psi:\tilde{X}'=\Spec_{Y'}(\mathscr{A}')&\longrightarrow\Spec_{Y'}(\mathscr{B})\subset\Spec_{Y'}(\psi^*\mathscr{A})\\
&=\tilde{X}(\D)\times_Y Y'\longrightarrow\tilde{X}(\D)
\end{split}
\]
is the composition of a finite morphism, a closed immersion and a base change of $\psi$, hence proper. On global sections, $r_{\D}\circ\Psi$ induces the identity of $A[Y,\D]$, so $r'=r_{\D}\circ\Psi$ is proper. Finally, the stalks of $\mathscr{A}$ and $\mathscr{A}'$ at the generic points are both $k(Y)[\omega^{\vee}\cap M]$, so $\Psi$, and hence $r'$, is birational.
\end{proof}

If $\psi$ is projective, then $\psi^*\D$ is a pp-divisor, $\tilde{X}'=\tilde{X}(\psi^*\D)$ and $r'=r_{\psi^*\D}$, by \cref{lemma: birational modification of a pp-div}\ref{lemma: birational modification of a pp-div part iv}.

\subsection{Proof of \cref{theorem: main theorem affine} for affine varieties.} The following result corresponds to the affine case of \cref{theorem: main theorem affine}.

\begin{proposition}\label{proposition: main theorem affine version}
    Every affine $T$-variety of complexity at most three admits an equivariant resolution of singularities.
\end{proposition}

\begin{proof}
Let $X$ be an affine $T$-variety of complexity at most three. Since $T$ is smooth, the action lift to the normalization and therefore we can assume that $X$ is normal. By \cite[Theorem 1.2(ii)]{MN25a}, $X\cong X(\D)$ for a pp-divisor $\D$ over a normal semi-projective variety $Y$ with $\dim(Y)\leq3$. Let $\psi:Y'\to Y$ and $B'$ be as in \cref{proposition: log resolution base} for $B:=\Supp(\D)$. Then, $\Supp(\psi^*\D)\leq B'$ by \cref{lemma: birational modification of a pp-div}\ref{lemma: birational modification of a pp-div part ii}. By \cref{proposition: proper pullback modification}, $r':\tilde{X}'\to X(\D)\cong X$ is a $T$-equivariant proper birational morphism from a normal $T$-variety, and $U_{B'}(\psi^*\D)\subset\tilde{X}'$ is a $T$-stable toroidal embedding without self-intersection. By \cref{proposition: Theorem 11* KKMSD73} there exists a resolution of singularities $\Psi:W\to\tilde{X}'$, and $r'\circ\Psi:W\to X$ is a resolution of singularities. By \cref{theorem: equivariant resolution toroidal}, the canonical resolution $f:\tilde{X}'^{\mathrm{can}}\to\tilde{X}'$ is $T$-equivariant and projective, so $r'\circ f$ is a $T$-equivariant resolution of singularities of $X$. When \cref{proposition: existence of snc pp-divisors} applies, $\tilde{X}'=\tilde{X}(\psi^*\D)$ with $\psi^*\D$ an SNC pp-divisor, and $r'=r_{\psi^*\D}$.
\end{proof}

\section{Equivariant resolution via toroidal embeddings: the general case}\label{Section: resolution general case}

In this section, we recall the theory of non-affine normal $T$-varieties and their description via \emph{divisorial fans}. We then construct a normal $T$-variety $\tilde{X}$ together with a proper birational $T$-equivariant morphism $r:\tilde{X}\to X$, which will play the role of $\tilde{X}(\D)$ in the affine case. Then, we introduce \emph{simple normal crossing divisorial fans}, and prove \cref{theorem: main theorem affine} in the general case.

\subsection{Divisorial fans and normal $T$-varieties.}\label{Section: divisorial fans and normal T-varieties}
To construct a divisorial fan, we allow empty coefficients in the definition of a pp-divisor. Following \cite[Section 2]{AHS08} and \cite[Section 3.1]{MN25b}, a pp-divisor over a normal semi-projective variety $Y$ is now allowed to be a formal sum
\[
\D:=\sum\Delta_D\otimes D
\]
where some $\Delta_D=\emptyset$ (with the rules $\emptyset+\Delta=\emptyset$ and $0\cdot\emptyset=\omega$), with the following conventions: $\bigcup_{\Delta_D=\emptyset}D$ is the support of an effective semiample divisor, the open subset $\Loc(\D):=Y\smallsetminus\bigcup_{\Delta_D=\emptyset}D$ is called the \emph{locus} of $\D$, and $\D$ is understood as the pp-divisor $\D|_{\Loc(\D)}:=\sum_{\Delta_D\neq\emptyset}\Delta_D\otimes(D\cap\Loc(\D))$ on $\Loc(\D)$, in the sense of the previous sections. Accordingly,
\[
\mathscr{A}(Y,\D):=\mathscr{A}(\Loc(\D),\D|_{\Loc(\D)}),\quad A[Y,\D]:=A[\Loc(\D),\D|_{\Loc(\D)}],
\]
$\tilde{X}(\D):=\Spec_{\Loc(\D)}(\mathscr{A}(Y,\D))$ and $X(\D):=\Spec(A[Y,\D])$. Besides, the respective morphisms $\pi_{\D}:\tilde{X}(\D)\to\Loc(\D)\subset Y$ and $r_{\D}:\tilde{X}(\D)\to X(\D)$ are as in \cref{Section: recall T-varieties}. More generally, given an open subvariety $V\subset Y$, we define the restriction of $\D$ to $V$ as $\D|_V:=\sum\Delta_D\otimes D|_V$. We say that $\D$ is an SNC pp-divisor if $\D|_{\Loc(\D)}$ is one, and put $U(\D):=U(\D|_{\Loc(\D)})$ and $Z(\D):=Z(\D|_{\Loc(\D)})$.

Let $\D$ and $\D'$ be two pp-divisors over the same normal semi-projective variety $Y$. Assume that $\Delta'_D\subset\Delta_D$ for every prime divisor $D$ on $Y$. Then, we have that $\Loc(\D')\subset\Loc(\D)$, the tail $\omega'$ of $\D'$ is contained in the tail $\omega$ of $\D$, and $\D(m)\leq\D'(m)$ on $\Loc(\D')$ for $m\in\omega^{\vee}\subset\omega'^{\vee}$. Hence $\mathscr{A}(Y,\D)|_{\Loc(\D')}\subset\mathscr{A}(Y,\D')$ and $A[Y,\D]\subset A[Y,\D']$ inside $k(Y)[M]$. These inclusions induce equivariant morphisms
\[
\iota:X(\D')\to X(\D)\quad\textrm{and}\quad
\tilde{\iota}:\tilde{X}(\D')\to\pi_{\D}^{-1}(\Loc(\D'))\subset\tilde{X}(\D),
\]
with $\pi_{\D}\circ\tilde{\iota}=\pi_{\D'}$ and $r_{\D}\circ\tilde{\iota}=\iota\circ r_{\D'}$ (both sides are morphisms to the affine scheme $X(\D)$ inducing the inclusion $A[Y,\D]\subset\H^0\left(\tilde{X}(\D'),\mathscr{O}_{\tilde{X}(\D')}\right)$). The pp-divisor $\D'$ is called a \emph{face} of $\D$, written $\D'\preceq\D$, if and only if $\Delta'_D\subset\Delta_D$ for every prime divisor $D$ on $Y$ and $\iota:X(\D')\to X(\D)$ is an open embedding. By \cite[Proposition 3.13]{MN25b}, this agrees with the combinatorial definition of \cite[Definition 5.1]{AHS08}.

A \emph{divisorial fan} over a normal semi-projective variety $Y$ is a finite set $\S$ of pp-divisors such that, for each pair $\D$ and $\D'$ in $\S$, the intersection
\[
\D\cap\D':=\sum(\Delta_D\cap\Delta'_D)\otimes D
\]
belongs to $\S$ and is a face of both $\D$ and $\D'$. The open embeddings $\iota_{\E\D}:X(\E)\to X(\D)$, $\E\preceq\D$ in $\S$, form a direct system: for $\mathfrak{F}\preceq\E\preceq\D$ one has $\iota_{\mathfrak{F}\D}=\iota_{\E\D}\circ\iota_{\mathfrak{F}\E}$, all maps being induced by inclusions of subalgebras of $k(Y)[M]$. The prevariety $X(\S):=\varinjlim_{\D\in\S}X(\D)$ is obtained by gluing the $X(\D)$ along the open subsets $X(\D\cap\D')$. Besides, it carries a $T$-action for which the $X(\D)\to X(\S)$ are equivariant open embeddings with $X(\D)\cap X(\D')=X(\D\cap\D')$ \cite[Theorem 5.3]{AHS08}, \cite[Section 3.4 and Proposition 3.20]{MN25b}. By \cite[Theorem 1.4]{MN25b}, every normal $T$-variety is equivariantly isomorphic to some $X(\S)$. For all divisorial fans considered below, we assume that $X(\S)$ is separated. The function field of $X(\S)$ is $k(Y)(M)$.

We define the \emph{support} of $\S$ as the reduced effective divisor
\[
\Supp(\S):=\sum D,
\]
where the sum runs over the set of prime divisors $D\text{ with }\Delta_D^{\D}\notin\{\emptyset,\omega_{\D}\}\text{ for some }\D\in\S$. Thus $\Supp(\D|_{\Loc(\D)})\leq\Supp(\S)|_{\Loc(\D)}$ for every $\D\in\S$.

\subsection{Partial desingularization.}\label{Section: partial desingularization non-affine}
Let $X$ be a normal $T$-variety. By \cite[Theorem 1.4]{MN25b}, there exists a divisorial fan $\S$ over a normal semi-projective variety $Y$ such that $X(\S)\cong X$ equivariantly. For every $\D\in\S$, the affinisation $r_{\D}:\tilde{X}(\D)\to X(\D)$ is a partial equivariant desingularization. A face relation $\E\preceq\D$ in $\S$ induces a $T$-equivariant open embedding $\iota_{\E\D}:X(\E)\to X(\D)$ and a $T$-equivariant morphism $\tilde{\iota}_{\E\D}:\tilde{X}(\E)\to\tilde{X}(\D)$. We prove that the latter is an open embedding using the following description of $\tilde{X}(\D)$.

\begin{lemma}\label{lemma: graph normalization}
Let $X,Y$ be varieties, $V\subset Y$ open, and $\tilde{X}$ a normal variety with a proper birational morphism $r:\tilde{X}\to X$ and an affine morphism $\pi:\tilde{X}\to V$. Let $\varpi:=\pi\circ r^{-1}:X\dashrightarrow Y$ and $\mathrm{Graph}_{\varpi}\subset X\times Y$ the closure of its graph. Then $g:=(r,\pi):\tilde{X}\to X\times Y$ is finite with image $\mathrm{Graph}_{\varpi}$, and $\tilde{X}\to\mathrm{Graph}_{\varpi}$ is the normalization of $\mathrm{Graph}_{\varpi}$. For every open $X_0\subset X$, $\mathrm{Graph}_{\varpi}\cap(X_0\times Y)=\mathrm{Graph}_{\varpi|_{X_0}}$.
\end{lemma}
\begin{proof}
The morphism $g$ is proper, since $r=\pr_1\circ g$ is proper and $\pr_1$ is separated \cite[Corollary II.4.8(e)]{Har77}. For affine open subsets $X_1\subset X$, $V_1\subset V$, $g^{-1}(X_1\times V_1)=\pi^{-1}(V_1)\times_X X_{1}$ is closed in the affine scheme $\pi^{-1}(V_{1})\times X_{1}$ ($X$ is separated), so $g:\tilde{X}\to X\times V$ is affine and proper, hence finite \cite[Exercise II.4.6]{Har77}. Its image $Z$ is closed in $X\times Y$, irreducible of dimension $\dim X$, and contains the graph of $\varpi$ over an open subset where $r$ is an isomorphism. Hence, it follows that $Z=\mathrm{Graph}_{\varpi}$. Being finite and birational from a normal variety, $\tilde{X}\to\mathrm{Graph}_{\varpi}$ is the normalization. The last assertion follows because the closure in the open subset $X_0\times Y$ is the intersection of the closure with it.
\end{proof}

We use two elementary facts about the normalization $\nu:Z^{\nu}\to Z$ of a variety: (a) every dominant morphism from a normal variety to $Z$ factors uniquely through $\nu$ \cite[Exercise II.3.8]{Har77}; and (b) if $\nu_i:Z_i\to Z$ ($i=1,2$) are finite birational morphisms from normal varieties and $\phi:Z_1\to Z_2$ satisfies $\nu_2\circ\phi=\nu_1$, then $\phi$ is an isomorphism (both $\nu_{i*}\mathscr{O}_{Z_i}$ are the integral closure of $\mathscr{O}_{Z}$ in $k(Z)$, and $\phi^{\#}$ is an $\mathscr{O}_{Z}$-algebra endomorphism of it, hence the identity).

\begin{lemma}\label{lemma: face open embedding}
Let $\E\preceq\D$ be a face relation of pp-divisors over $Y$. Then $\tilde{\iota}_{\E\D}:\tilde{X}(\E)\to\tilde{X}(\D)$ is a $T$-equivariant open embedding with image $r_{\D}^{-1}(X(\E))$. More generally, let $\psi:Y'\to Y$ be a proper birational morphism from a smooth variety $Y'$, $B'$ a simple normal crossing divisor on $Y'$, $\Loc'(\D):=\psi^{-1}(\Loc(\D))$, $\D':=(\psi|_{\Loc'(\D)})^*(\D|_{\Loc(\D)})$, and similarly $\E'$, and assume $\Supp(\D'),\Supp(\E')\leq B'$. Then the morphism $\tilde{\iota}'_{\E\D}:\tilde{X}'(\E)\to\tilde{X}'(\D)$ between the varieties of \cref{proposition: proper pullback modification}, induced by $\mathscr{A}(\Loc'(\D),\D')|_{\Loc'(\E)}\subset\mathscr{A}(\Loc'(\E),\E')$, is a $T$-equivariant open embedding with image $(r'_{\D})^{-1}(X(\E))$.
\end{lemma}
\begin{proof}
The first assertion is the case $\psi=\mathrm{id}_Y$ (then $Y$ need not be smooth and only \cref{lemma: graph normalization} is used). The inclusion of algebras exists because $\E(m)\geq\D(m)$ on $\Loc(\E)$ for $m\in\omega_{\D}^{\vee}$ and pull-back preserves effectivity. Identify $X(\E)$ with the open subset $\iota_{\E\D}(X(\E))\subset X(\D)$. The rational maps $\varpi'_{\D}:=\pi'_{\D}\circ(r'_{\D})^{-1}:X(\D)\dashrightarrow Y'$ and $\varpi'_{\E}$ both correspond to $k(Y')=k(Y)\subset k(Y)(M)$, so $\varpi'_{\E}=\varpi'_{\D}|_{X(\E)}$. By \cref{lemma: graph normalization}, applied to $(\tilde{X}'(\D),r'_{\D},\pi'_{\D})$, the open subset $(r'_{\D})^{-1}(X(\E))=g_{\D}^{-1}(X(\E)\times Y')$ is the normalization of $\mathrm{Graph}_{\varpi'_{\D}}\cap(X(\E)\times Y')=\mathrm{Graph}_{\varpi'_{\E}}$, and by the same lemma applied to $(\tilde{X}'(\E),r'_{\E},\pi'_{\E})$, so is $\tilde{X}'(\E)$. The relations $\pi'_{\D}\circ\tilde{\iota}'_{\E\D}=\pi'_{\E}$ and $r'_{\D}\circ\tilde{\iota}'_{\E\D}=\iota_{\E\D}\circ r'_{\E}$ (the latter because both are morphisms to the affine scheme $X(\D)$ inducing $A[Y,\D]\subset A[Y,\E]$ on global sections) give $g_{\D}\circ\tilde{\iota}'_{\E\D}=g_{\E}$. Thus $\tilde{\iota}'_{\E\D}$ maps $\tilde{X}'(\E)$ into $(r'_{\D})^{-1}(X(\E))$ compatibly with the normalization morphisms to $\mathrm{Graph}_{\varpi'_{\E}}$; and, therefore by fact (b), it is an isomorphism onto this open subset. The morphism is equivariant since the inclusion of algebras is graded.
\end{proof}

For any $\mathfrak{F}\preceq\E\preceq\D$ in $\S$, we have $\iota_{\mathfrak{F}\D}=\iota_{\E\D}\circ\iota_{\mathfrak{F}\E}$ and $\tilde{\iota}_{\mathfrak{F}\D}=\tilde{\iota}_{\E\D}\circ\tilde{\iota}_{\mathfrak{F}\E}$, all these maps being induced by inclusions of subalgebras of $k(Y)[M]$. They are compatible with the morphisms $\pi$ and $r$: $\pi_{\D}\circ\tilde{\iota}_{\E\D}=\pi_{\E}$ and $r_{\D}\circ\tilde{\iota}_{\E\D}=\iota_{\E\D}\circ r_{\E}$ (\cref{Section: divisorial fans and normal T-varieties}).

Thus, the face relations in $\S$ give the normal $T$-varieties $\tilde{X}(\D)$ the structure of a direct system, and the equivariant open embeddings $\tilde{\iota}_{\E\D}$ give gluing data, as for $X(\D)$ \cite[Section 3.4]{MN25b}. By \cref{lemma: face open embedding}, inside $\tilde{X}(\D)$ we have
\[
\begin{split}
\tilde{X}(\D\cap\D')\cap\tilde{X}(\D\cap\D'')
&=r_{\D}^{-1}(X(\D\cap\D'\cap\D''))\\
&=\tilde{X}(\D\cap\D'\cap\D''),
\end{split}
\]
so the cocycle condition holds \cite[Exercise II.2.12]{Har77}. Moreover, the $T$-equivariant proper birational maps $r_{\D}:\tilde{X}(\D)\to X(\D)$ glue to a $T$-equivariant morphism
\[
r_{\S}:\tilde{X}(\S):=\varinjlim_{\D\in\S}\tilde{X}(\D)\longrightarrow\varinjlim_{\D\in\S}X(\D)=X(\S)
\]
of $T$-prevarieties. Likewise, the $T$-invariant morphisms $\pi_{\D}:\tilde{X}(\D)\to\Loc(\D)\subset Y$ glue to a $T$-invariant morphism $\pi_{\S}:\tilde{X}(\S)\to Y$. Therefore, the partial equivariant desingularizations $r_{\D}:\tilde{X}(\D)\to X(\D)$ glue to a partial equivariant desingularization $r_{\S}:\tilde{X}(\S)\to X(\S)$.

\begin{proposition}\label{proposition: partial equivariant desingularization}
Let $\S$ be a divisorial fan over a normal semi-projective variety $Y$. Then $\tilde{X}(\S)$ is a normal $T$-variety. Moreover, the morphism $r_{\S}:\tilde{X}(\S)\to X(\S)$ is $T$-equivariant, proper and birational.
\end{proposition}
\begin{proof}
First, observe that for any pair $\D$ and $\E$ in $\S$, by \cref{lemma: face open embedding},
\[
r_{\D}^{-1}(X(\E))=r_{\D}^{-1}(X(\D)\cap X(\E))=r_{\D}^{-1}(X(\D\cap\E))=\tilde{X}(\D\cap\E).
\]
Consequently, $r_{\S}^{-1}(X(\D))=\tilde{X}(\D)$ for every $\D\in\S$. Since each $r_{\D}$ is proper, \cite[Corollary II.4.8(f)]{Har77} implies that $r_{\S}$ is proper. As $X(\S)$ is separated, so is $\tilde{X}(\S)$. It is of finite type, normal, and irreducible, being the finite union of the irreducible open subsets $\tilde{X}(\D)$, any two of which meet in $\tilde{X}(\D\cap\D')\neq\emptyset$ (since $X(\D)\cap X(\D')=X(\D\cap\D')$ is a non-empty open subset of the irreducible $X(\S)$). The $T$-actions glue because the gluing maps are equivariant, and the action is effective since it is effective on each $\tilde{X}(\D)$. Hence $\tilde{X}(\S)$ is a normal $T$-variety.

The morphism $r_{\S}$ is equivariant because it is obtained by gluing $T$-equivariant morphisms, and it is birational because its restriction over the dense open subset $X(\D)$ of $X(\S)$ is the birational morphism $r_{\D}$.
\end{proof}

\subsection{Simple normal crossing divisorial fans.}
As in the affine case, we can modify a divisorial fan to obtain a partial equivariant desingularization which admits a toroidal embedding without self-intersection.
\begin{definition}\label{definition: snc divisorial fan}
A divisorial fan $\S$ over a normal semi-projective variety $Y$ is called \emph{simple normal crossing} if $Y$ is smooth and $\Supp(\S)$ is a simple normal crossing divisor.
\end{definition}
\begin{proposition}\label{proposition: snc divisorial fans iff snc pp-div}
If $\S$ is a simple normal crossing divisorial fan, then every $\D\in\S$ is a simple normal crossing pp-divisor. More precisely, for every $\D\in\S$ the divisor $B_{\D}:=\Supp(\S)|_{\Loc(\D)}$ is a simple normal crossing divisor on the smooth variety $\Loc(\D)$ with $\Supp(\D|_{\Loc(\D)})\leq B_{\D}$.
\end{proposition}
\begin{proof}
Since $Y$ is smooth, so is its open subset $\Loc(\D)$, and the components of $\Supp(\D|_{\Loc(\D)})$ are among those of the simple normal crossing divisor $B_{\D}$, which is the restriction of $\Supp(\S)$ to an open subset. Hence both are simple normal crossing divisors.
\end{proof}

Let $\S$ be a simple normal crossing divisorial fan. For $\D\in\S$, consider the $T$-stable open subset $U_{B_{\D}}(\D|_{\Loc(\D)})\subset\tilde{X}(\D)$ of \cref{Section: snc ppdivisors}, whose closed points are precisely the points $x$ with $\pi_{\D}(x)\notin\Supp(\S)$ and $\dim(T\cdot x)=\dim(T)$. For $\E\preceq\D$ in $\S$, the open embedding $\tilde{\iota}_{\E\D}$ is $T$-equivariant and compatible with the projections to $Y$, so this description shows that $U_{B_{\E}}(\E|_{\Loc(\E)})=U_{B_{\D}}(\D|_{\Loc(\D)})\cap\tilde{X}(\E)$. Hence these open subsets glue to a $T$-stable open subset $U(\S)\subset\tilde{X}(\S)$, and we put $Z(\S):=\tilde{X}(\S)\smallsetminus U(\S)$, so that, for every $\D\in\S$,
\[
\begin{split}
Z(\S)\cap\tilde{X}(\D)&=Z_{B_{\D}}(\D|_{\Loc(\D)})\\
&=\pi_{\D}^{-1}(\Supp(\S)\cap\Loc(\D))\cup\bigcup_{\varrho\in\omega_{\D}(1)}\tilde{E}_{\varrho}.
\end{split}
\]
Note that $Z(\S)\cap\tilde{X}(\D)$ contains $Z(\D)$, but may be strictly larger: a prime divisor of $Y$ may have trivial coefficient in $\D$ and a non-trivial one in another $\D'\in\S$. This causes no difficulty, since \cref{corollary: snc pp div and toroidal embeddings} allows an arbitrary simple normal crossing divisor $B\geq\Supp(\D)$.

\begin{proposition}\label{proposition: snc div fan and toroidal embeddings}
If $\S$ is a simple normal crossing divisorial fan, then $U(\S)\subset\tilde{X}(\S)$ is a toroidal embedding without self-intersection.
\end{proposition}
\begin{proof}
By \cref{corollary: snc pp div and toroidal embeddings,proposition: snc divisorial fans iff snc pp-div}, applied to $\D|_{\Loc(\D)}$ and $B_{\D}$, each inclusion $U(\S)\cap\tilde{X}(\D)=U_{B_{\D}}(\D|_{\Loc(\D)})\subset\tilde{X}(\D)$, with $\D\in\S$, is a toroidal embedding without self-intersection. In particular, $U(\S)$ is smooth, dense (each $U_{B_{\D}}(\D|_{\Loc(\D)})$ is dense in $\tilde{X}(\D)$) and $T$-stable (each $U_{B_{\D}}(\D|_{\Loc(\D)})$ is). Since the $T$-stable open subvarieties $\tilde{X}(\D)$ cover $\tilde{X}(\S)$, every closed point $x\in\tilde{X}(\S)$ lies in some $\tilde{X}(\D)$. Consequently,
\[
\hat{\mathscr{O}}_{\tilde{X}(\S),x}\cong\hat{\mathscr{O}}_{\tilde{X}(\D),x}\cong\hat{\mathscr{O}}_{X_{\tau},t},
\]
where $X_{\tau}$ is an affine $T'$-toric variety.

Moreover, since $Z(\S)\cap\tilde{X}(\D)=Z_{B_{\D}}(\D|_{\Loc(\D)})$, the ideals of $\tilde{X}(\S)\smallsetminus U(\S)$ and $\tilde{X}(\D)\smallsetminus U_{B_{\D}}(\D|_{\Loc(\D)})$ generate the same ideal in $\hat{\mathscr{O}}_{X_{\tau},t}$. Thus, the local toroidal structure on $\tilde{X}(\D)$ also identifies the boundary of $U(\S)$ with the toric boundary, and the étale charts of $\tilde{X}(\D)$ are étale charts of $\tilde{X}(\S)$. Finally, if $E$ is an irreducible component of $Z(\S)$ and $x\in E$ is a closed point with $x\in\tilde{X}(\D)$, then $E\cap\tilde{X}(\D)$ is an irreducible component of $Z_{B_{\D}}(\D|_{\Loc(\D)})$, which is normal by \cref{proposition: snc pp div and toroidal embeddings}. Hence $\mathscr{O}_{E,x}$ is normal, and the irreducible components of $Z(\S)$ are normal also by \cref{proposition: snc pp div and toroidal embeddings}.
\end{proof}

\subsection{Proof of \cref{theorem: main theorem affine}.}
By \cite[Theorem 1.4]{MN25b}, every normal $T$-variety arises from a divisorial fan over a normal semi-projective variety. As in the affine case, this algebro-geometric combinatorial description is not unique: two non-isomorphic divisorial fans $\S$ and $\S'$ can give rise to isomorphic normal $T$-varieties, as the following result shows.

\begin{lemma}\label{lemma: birational modification of a div fan}
Let $\S$ be a divisorial fan over a normal semi-projective variety $Y$ and $\psi:Y'\to Y$ a projective birational morphism. For $\D\in\S$ let $\psi^*\D$ be the polyhedral divisor on $Y'$ with empty coefficients along the prime divisors contained in $Y'\smallsetminus\Loc'(\D)$, $\Loc'(\D):=\psi^{-1}(\Loc(\D))$, and equal to $(\psi|_{\Loc'(\D)})^*(\D|_{\Loc(\D)})$ on $\Loc'(\D)$. Then
\[
\psi^*\S:=\{\psi^*\D\mid\D\in\S\}
\]
is a divisorial fan over $Y'$ with $\psi^*(\D\cap\D')=\psi^*\D\cap\psi^*\D'$, and the normal $T$-varieties $X(\psi^*\S)$ and $X(\S)$ are isomorphic as $T$-varieties.
\end{lemma}
\begin{proof}
\emph{Step 1: $\psi^*\D$ is a pp-divisor.} By \cref{lemma: birational modification of a pp-div}\ref{lemma: birational modification of a pp-div part iv}, $Y'$ is semi-projective. If $H$ is an effective semiample divisor on $Y$ with support $Y\smallsetminus\Loc(\D)$, then $\psi^*H$ is effective and semiample with support $Y'\smallsetminus\Loc'(\D)$, so the convention on empty coefficients is satisfied, with $\Loc(\psi^*\D)=\Loc'(\D)$. Since $\psi|_{\Loc'(\D)}:\Loc'(\D)\to\Loc(\D)$ is projective and birational, \cref{lemma: birational modification of a pp-div}\ref{lemma: birational modification of a pp-div part iv} shows that $\psi^*\D|_{\Loc'(\D)}$ is a pp-divisor and $X(\psi^*\D)=X(\D)$, both being the spectrum of the same subalgebra $A[Y,\D]$ of $k(Y)[M]$.

\emph{Step 2: a valuative description of the coefficients.} Let $F$ be a prime divisor of $Y'$ meeting $\Loc'(\D)$, with valuation $\nu_F$ on $k(Y)$, and let $v\in N_{\Q}$ and $l>0$ with $lv\in N$. The map
\[
\nu:k(Y)[M]\smallsetminus\{0\}\to\Z,\qquad
\nu\left(\sum_m f_m\chi^m\right):=\min_m\bigl(l\nu_F(f_m)+\langle m,lv\rangle\bigr),
\]
is a discrete valuation of $k(Y)(M)$, the monomial valuation attached to $l\nu_F$ and the weight $lv\in N$: $\nu(f+g)\geq\min(\nu(f),\nu(g))$ and $\nu(fg)\geq\nu(f)+\nu(g)$ are clear, and if $m_1,m_2$ are the largest indices (for a total group ordering of $M$) where the minima defining $\nu(f),\nu(g)$ are attained, then the coefficient of $\chi^{m_1+m_2}$ in $fg$ is $f_{m_1}g_{m_2}$ plus terms of strictly larger $\nu_F$-value, so $\nu(fg)\leq\nu(f)+\nu(g)$. We claim that
\[
v\in\Delta_F^{\psi^*\D}\quad\Longleftrightarrow\quad\nu\geq0\text{ on }A[Y,\D].
\]
Write $h_F:=h_{\Delta_F^{\psi^*\D}}$, so that $h_F(m)$ is the coefficient of $F$ in $\psi^*\D(m)$, and recall that $A[Y,\D]=A[\Loc'(\D),\psi^*\D]$ by Step 1. If $v\in\Delta_F^{\psi^*\D}$ and $f\chi^m\in A[Y,\D]$, then $\nu_F(f)+\lfloor h_F(m)\rfloor\geq0$ and $\langle m,v\rangle\geq h_F(m)$, so $l\nu_F(f)+\langle m,lv\rangle\geq0$. Since $\nu$ is a minimum over the homogeneous components, $\nu\geq0$ on $A[Y,\D]$. Conversely, assume $\nu\geq0$ on $A[Y,\D]$ and let $m\in\omega_{\D}^{\vee}\cap M$. Choose $n>0$ such that $L:=n\psi^*\D(m)$ is Cartier on $\Loc'(\D)$. As $L$ is semiample, there is $f\in\H^0(\Loc'(\D),\mathscr{O}(L))$ with $F\not\subset\Supp(\div(f)+L)$, i.e. $\nu_F(f)=-nh_F(m)$. Then $f\chi^{nm}\in A[Y,\D]$, so $0\leq\nu(f\chi^{nm})=l(-nh_F(m)+\langle nm,v\rangle)$, i.e. $\langle m,v\rangle\geq h_F(m)$. This holds for all $m\in\omega_{\D}^{\vee}\cap M$, hence for all $m\in\omega_{\D}^{\vee}$ by homogeneity, and $v\in\Delta_F^{\psi^*\D}$ since a polyhedron is determined by its support function.

\emph{Step 3: intersections.} Let $\D,\D'\in\S$ and $\E:=\D\cap\D'\in\S$, a face of both. Note that $\Loc(\E)\subset\Loc(\D)\cap\Loc(\D')$, possibly strictly (an intersection of non-empty polyhedra may be empty). Let $F$ be a prime divisor of $Y'$.

(a) If $F$ misses $\Loc'(\D)$ or $\Loc'(\D')$, then $F$ misses $\Loc'(\E)$ and both $\psi^*\E$ and $\psi^*\D\cap\psi^*\D'$ have coefficient $\emptyset$ along $F$.

(b) Assume $F$ meets $\Loc'(\D)$ and $\Loc'(\D')$, let $v\in N_{\Q}$ and $\nu$ be as in Step 2, and suppose $\nu\geq0$ on $A[Y,\D]$ and on $A[Y,\D']$. Let $R_{\nu}\subset K:=k(Y)(M)$ be the valuation ring of $\nu$. The inclusions $A[Y,\D]\subset R_{\nu}$ and $A[Y,\D']\subset R_{\nu}$ define morphisms $\Spec(R_{\nu})\to X(\D)\subset X(\S)$ and $\Spec(R_{\nu})\to X(\D')\subset X(\S)$ which agree on $\Spec(K)$. Since $X(\S)$ is separated, they coincide \cite[Theorem II.4.3]{Har77}. Hence the image of the closed point lies in $X(\D)\cap X(\D')=X(\E)$, the common morphism factors through the open subset $X(\E)$, and the corresponding ring homomorphism $A[Y,\E]=\mathscr{O}(X(\E))\to R_{\nu}$ is the inclusion $A[Y,\E]\subset K$. Thus, $\nu\geq0$ on $A[Y,\E]$.

(c) If moreover $F$ misses $\Loc'(\E)$, this is impossible: let $H$ be an effective semiample divisor on $Y$ with support $Y\smallsetminus\Loc(\E)$ and $n>0$ with $nH$ Cartier and the $Y_h$, $h\in\H^0(Y,\mathscr{O}_Y(nH))$, covering $Y$. Choose $h$ with $\psi(F)\not\subset\Supp(\div(h)+nH)$, which does not vanish at the generic point of $\psi(F)$. Then $h\in\mathscr{O}(\Loc(\E))=A[Y,\E]_0$, so $\nu_F(h)\geq0$. But $F\subset\psi^{-1}(\Supp(H))=\Supp(\psi^*H)$, so the coefficient of $F$ in $\psi^*(nH)$ is positive and that of $F$ in $\psi^*(\div(h)+nH)=\div(h)+\psi^*(nH)$ is zero, whence $\nu_F(h)<0$. Hence $\Delta_F^{\psi^*\D}\cap\Delta_F^{\psi^*\D'}=\emptyset$ is the coefficient of $\psi^*\E$ along $F$.

(d) If $F$ meets $\Loc'(\E)$, then (b) gives $v\in\Delta_F^{\psi^*\E}$ whenever $v\in\Delta_F^{\psi^*\D}\cap\Delta_F^{\psi^*\D'}$, and conversely $v\in\Delta_F^{\psi^*\E}$ gives $\nu\geq0$ on $A[Y,\E]\supset A[Y,\D],A[Y,\D']$. This proves $\psi^*(\D\cap\D')=\psi^*\D\cap\psi^*\D'$.

\emph{Step 4: faces.} For $\E\preceq\D$ in $\S$, Step 2 and $A[Y,\D]\subset A[Y,\E]$ give $\Delta_F^{\psi^*\E}\subset\Delta_F^{\psi^*\D}$ for all $F$, and the morphism $X(\psi^*\E)\to X(\psi^*\D)$ induced by the inclusion of algebras is, under the identifications $X(\psi^*\E)=X(\E)$ and $X(\psi^*\D)=X(\D)$ of Step 1, the open embedding $\iota_{\E\D}$. Hence $\psi^*\S$ is a divisorial fan over $Y'$, and since the $X(\D)$ cover $X(\S)$ and the identifications $X(\psi^*\D)=X(\D)$ are compatible with the gluing maps, $X(\S)\cong X(\psi^*\S)$ equivariantly.
\end{proof}

Let $\S$ be a divisorial fan over a normal semi-projective variety $Y$ that gives rise to $X$. By \cref{lemma: birational modification of a div fan}, any projective birational morphism $\psi:Y'\to Y$ induces a $T$-equivariant isomorphism $X(\S')\cong X(\S)$, where $\S':=\psi^*\S$. However, the associated $T$-varieties $\tilde{X}(\S)$ and $\tilde{X}(\S')$ need not be isomorphic, and it may happen that only one admits a toroidal embedding, for instance when its defining divisorial fan is simple normal crossing. The existence of such a divisorial fan requires a smooth base $Y$ and is therefore not guaranteed in arbitrary complexity, although it can be ensured in low complexity under the following assumptions.

\begin{proposition}\label{proposition: existence of snc divisorial fan}
Let $X$ be a normal $T$-variety of complexity at most three, and let $\S'$ be a divisorial fan over a normal semi-projective variety $Y$ such that $X\cong X(\S)$. If $\dim(Y)\leq2$, or if the morphism $\psi$ of \cref{proposition: log resolution base} for $(Y,\Supp(\S))$ can be chosen projective, then there exists a simple normal crossing divisorial fan $\S'$ such that $X\cong X(\S')$.
\end{proposition}
\begin{proof}
For $\D\in\S$, the open subset $X(\D)\subset X$ has complexity $\dim\Loc(\D)=\dim(Y)$, so $\dim(Y)\leq3$. Let $\psi:Y'\to Y$ and $B$ be as in \cref{proposition: log resolution base} for $(Y,\Supp(\S))$, with $\psi$ projective. By \cref{lemma: birational modification of a div fan}, $\S':=\psi^*\S$ is a divisorial fan over the smooth variety $Y'$ with $X(\S)\cong X(\S')\cong X$. A prime divisor $F$ of $Y'$ with $\Delta_F^{\psi^*\D}\notin\{\emptyset,\omega_{\D}\}$ meets $\Loc'(\D)$, and by \cref{lemma: birational modification of a pp-div}\ref{lemma: birational modification of a pp-div part ii}, applied to $\psi|_{\Loc'(\D)}$, it is either exceptional or the strict transform of a component of $\Supp(\D|_{\Loc(\D)})\leq\Supp(\S)$. In both cases $F\subset\psi^{-1}(\Supp(\S)\cup\Sing(Y))\subset B$, since the set $W$ of \cref{Section: pullbacks} for $\psi|_{\Loc'(\D)}$ is contained in $\Supp(\S)\cup\Sing(Y)$. Hence the components of $\Supp(\S')$ are among those of $B$, and $\S'$ is a simple normal crossing divisorial fan.
\end{proof}

The following analogue of \cref{proposition: proper pullback modification} removes the projectivity assumption and shows that \cref{theorem: main theorem affine} holds without a projectivity assumption.

\begin{proposition}\label{proposition: proper fan modification}
Let $\S$ be a divisorial fan over a normal semi-projective variety $Y$, $\psi:Y'\to Y$ a proper birational morphism from a smooth variety $Y'$, and $B'$ a simple normal crossing divisor on $Y'$ such that, in the notation of \cref{lemma: face open embedding}, $\Supp(\D')\leq B'|_{\Loc'(\D)}$ for every $\D\in\S$. Then the varieties $\tilde{X}'(\D)$, $\D\in\S$, glue along the open embeddings $\tilde{\iota}'_{\E\D}$ to a normal $T$-variety $\tilde{X}_{\psi}(\S)$ endowed with a $T$-equivariant proper birational morphism $r'_{\S}:\tilde{X}_{\psi}(\S)\to X(\S)$ with $(r'_{\S})^{-1}(X(\D))=\tilde{X}'(\D)$, and the open subsets $U_{B'}(\D')\subset\tilde{X}'(\D)$ glue to a $T$-stable open subset $U_{\psi}(\S)\subset\tilde{X}_{\psi}(\S)$ which is a toroidal embedding without self-intersection.
\end{proposition}
\begin{proof}
By \cref{proposition: proper pullback modification}, applied to $\D|_{\Loc(\D)}$, $\psi|_{\Loc'(\D)}$ and $B'|_{\Loc'(\D)}$, each $\tilde{X}'(\D)$ is a normal $T$-variety with a $T$-equivariant proper birational morphism $r'_{\D}:\tilde{X}'(\D)\to X(\D)$, and $U_{B'}(\D')\subset\tilde{X}'(\D)$ is a toroidal embedding without self-intersection. By \cref{lemma: face open embedding}, $\tilde{X}'(\D\cap\D')=(r'_{\D})^{-1}(X(\D\cap\D'))$ is open in $\tilde{X}'(\D)$, and inside $\tilde{X}'(\D)$ we have $\tilde{X}'(\D\cap\D')\cap\tilde{X}'(\D\cap\D'')=\tilde{X}'(\D\cap\D'\cap\D'')$. Thus the gluing is well defined, and the $r'_{\D}$ glue to $r'_{\S}$ with $(r'_{\S})^{-1}(X(\D))=\tilde{X}'(\D)$. The proof of \cref{proposition: partial equivariant desingularization} now applies with $r'_{\D}$ in place of $r_{\D}$, and shows that $\tilde{X}_{\psi}(\S)$ is a normal $T$-variety and $r'_{\S}$ is $T$-equivariant, proper and birational. As in Section 4.3, the closed points of $U_{B'}(\D')$ are those $x$ with $\pi'_{\D}(x)\notin B'$ and $\dim(T\cdot x)=\dim(T)$, and the $\tilde{\iota}'_{\E\D}$ are equivariant and over $Y'$. Hence $U_{B'}(\E')=U_{B'}(\D')\cap\tilde{X}'(\E)$, the $U_{B'}(\D')$ glue to $U_{\psi}(\S)$, and the proof of \cref{proposition: snc div fan and toroidal embeddings} applies with $Z_{B'}(\D')$ in place of $Z_{B_{\D}}(\D|_{\Loc(\D)})$.
\end{proof}

\begin{theoremff}{theorem: main theorem affine}
Every $T$-variety of complexity at most three admits an equivariant resolution of singularities.
\end{theoremff}
\begin{proof}
Let $X$ be a $T$-variety of complexity at most three. As in the proof of \cref{proposition: main theorem affine version}, the action lift to the normalization, since $T$ is smooth, and therefore we can assume that $X$ is normal.

Let $X$ be a normal $T$-variety of complexity at most three. By \cite[Theorem 1.4]{MN25b}, there exists a divisorial fan $\S$ over a normal semi-projective variety $Y$ such that $X\cong X(\S)$ with $\dim(Y)\leq3$ (as in the proof of \cref{proposition: existence of snc divisorial fan}). Let $\psi:Y'\to Y$ and $B'$ be as in \cref{proposition: log resolution base} for $(Y,\Supp(\S))$. For $\D\in\S$, the argument of \cref{proposition: existence of snc divisorial fan} shows that every component of $\Supp(\D')$ lies in $\psi^{-1}(\Supp(\S)\cup\Sing(Y))\cap\Loc'(\D)\subset B'|_{\Loc'(\D)}$, so the hypotheses of \cref{proposition: proper fan modification} hold. Hence $r'_{\S}:\tilde{X}_{\psi}(\S)\to X(\S)\cong X$ is a $T$-equivariant proper birational morphism from a normal $T$-variety, and the $T$-stable open subvariety $U_{\psi}(\S)\subset\tilde{X}_{\psi}(\S)$ is a toroidal embedding without self-intersection. By \cref{proposition: Theorem 11* KKMSD73}, there exists a resolution of singularities $\Psi:W\to\tilde{X}_{\psi}(\S)$, and the composition $r'_{\S}\circ\Psi:W\to X$ is a desingularization. By \cref{theorem: equivariant resolution toroidal}, the canonical resolution $f$ of $\tilde{X}_{\psi}(\S)$ is $T$-equivariant and projective, so $r'_{\S}\circ f$ is a $T$-equivariant resolution of singularities of $X$. When \cref{proposition: existence of snc divisorial fan} applies, $\psi^*\S$ is a simple normal crossing divisorial fan, $\tilde{X}_{\psi}(\S)=\tilde{X}(\psi^*\S)$ and $U_{\psi}(\S)=U(\psi^*\S)$ by \cref{proposition: snc div fan and toroidal embeddings}. This proves \cref{theorem: main theorem affine}. 
\end{proof}

\begin{remark}
We do not know whether the resulting resolutions are projective over $X$. We know that $\Psi$ and $f$ are projective, but $r'_{\S}$ (and $r'$) is only shown to be proper, as is the log resolution of \cref{proposition: log resolution base} in dimension three.
\end{remark}

\section{The divisorial fan of the resolution}\label{Section: combinatorial resolution}
In this section, we give a combinatorial description of the resolutions and we finish by providing a concrete example.

Before going on with the section, we recall that a \emph{fan} is a finite set of pointed cones closed under taking faces, and any two of them meet in a common face. A fan $\Sigma'$ \emph{refines} a fan $\Sigma$ if both have the same support and every cone of $\Sigma'$ lies in a cone of $\Sigma$. A cone is \emph{regular} if it is generated by part of a basis of the lattice.

\subsection{Notations and settings.} Throughout, $Y$ is a smooth variety of dimension $n$, $\D=\sum\Delta_D\otimes D$ a polyhedral divisor on $Y$ with pointed tail $\omega$, and $B=\sum_{i\in I}D_i$ a simple normal crossing divisor with $\Supp(\D)\leq B$, as in \cref{proposition: snc pp div and toroidal embeddings}. We write $\Delta_i:=\Delta_{D_i}$, $\tilde{X}:=\Spec_Y(\mathscr{A}(Y,\D))$, $D_J:=\bigcap_{j\in J}D_j$ for $J\subset I$, $\tilde{N}:=N\times\Z^I$ with standard basis $(e_i)$ of $\Z^I$, $p:\tilde{N}_{\Q}\to\Q^I$ for the projection, and regard $N_{\Q}\times\Q^J$ as a subspace of $\tilde{N}_{\Q}$. For a cone $\tau\subset N_{\Q}\times\Q^I_{\geq0}$, we put $J(\tau):=\{i\mid\lambda_i\not\equiv0\text{ on }\tau\}$, where $\lambda_i$ is the $i$-th coordinate, and $\tau_0:=\tau\cap(N_{\Q}\times\{0\})$.

\subsection{The Cayley fan and Cayley subdivisions.}
For $J\subset I$ such that $D_J\neq\emptyset$, the \emph{Cayley cone}
\[
\sigma_J:=\operatorname{cone}\left(\omega\times\{0\}\cup\bigcup_{j\in J}\Delta_j\times\{e_j\}\right)\subset N_{\Q}\times\Q^J
\]
is the cone $\sigma$ of \cref{proposition: snc pp div and toroidal embeddings}\ref{proposition: snc pp div and toroidal embeddings part ii} at any point of $D_J$ lying on no other component of $B$. By Step 1 of that proof, its elements are the $(w+\sum_j\lambda_jv_j,\lambda)$ with $w\in\omega$, $\lambda\in\Q^J_{\geq0}$, $v_j\in\Delta_j$. Hence its slice at $\lambda$ is $\sigma_J\cap(N_{\Q}\times\{\lambda\})=(\sum_j\lambda_j\Delta_j+\omega)\times\{\lambda\}$, and $\sigma_J\cap(N_{\Q}\times\Q^{J'})=\sigma_{J'}$ is the face of $\sigma_J$ cut out by $(0,\sum_{j\in J\smallsetminus J'}e_j)$ for $J'\subset J$. A face of $\sigma_J$ cut out by $(m,a)\in\sigma_J^{\vee}$ (so $m\in\omega^{\vee}$, $a_j\geq-h_{\Delta_j}(m)$) consists of the elements $(w+\sum_j\lambda_jv_j,\lambda)$ on which the functional, a sum of the non-negative terms $\langle m,w\rangle$ and $\lambda_j(\langle m,v_j\rangle+a_j)$, vanishes, i.e. $w\in\omega_{m}:=\omega\cap m^{\perp}$ and, for each $j$ with $\lambda_j>0$, $a_j=-h_{\Delta_j}(m)$ and $v_j\in\Delta_{j,m}:=\{v\in\Delta_j\mid\langle m,v\rangle=h_{\Delta_j}(m)\}$. Thus every face is of the form
\[
\mathrm{cone}\left(\omega_{m}\times\{0\}\cup\bigcup_{j\in J'}\Delta_{j,m}\times\{e_j\}\right)
\]
 for some $J'\subset J$, the Cayley cone of the faces $\omega_{m},\Delta_{j,m}$. Consequently, the faces of the $\sigma_J$ are finitely many, $\sigma_J\cap\sigma_{J'}=\sigma_{J\cap J'}$ (both sides equal $\sigma_J\cap\sigma_{J'}\cap(N_{\Q}\times\Q^{J\cap J'})$), and the intersection of a face of $\sigma_J$ with a face of $\sigma_{J'}$ is an intersection of two faces of $\sigma_{J\cap J'}$, hence a common face. Thus, the set $\Sigma(\D,B)$ of all faces of all $\sigma_J$, $D_J\neq\emptyset$, is a fan, the \emph{Cayley fan} of $(\D,B)$, and a \emph{subdivision} of $(\D,B)$ is a fan $\Sigma$ refining it. For $\tau\in\Sigma$ one has $\tau\subset\sigma_{J(\tau)}$ and $D_{J(\tau)}\neq\emptyset$ (if $\tau\subset\sigma_J$ then $J(\tau)\subset J$ and $\tau\subset\sigma_J\cap(N_{\Q}\times\Q^{J(\tau)})$). For $N=0$ (the trivial polyhedral divisor $\D_0$ on $Y$ with tail $\{0\}$), the Cayley fan is the \emph{orthant fan} $\Gamma_B:=\{\text{faces of }\Q^J_{\geq0}\mid D_J\neq\emptyset\}$ in $\Q^I$.

\begin{definition}\label{definition: Cayley subdivision}
Let $\Gamma$ be a subdivision of $\Gamma_B$, with primitive ray generators $n_{\rho}\in\Z^I_{\geq0}$, $\rho\in\Gamma(1)$. A subdivision $\Sigma$ of $\Sigma(\D,B)$ is \emph{Cayley over $\Gamma$} if, for every $\tau\in\Sigma$, the cone $p(\tau)$ belongs to $\Gamma$ and
\[
\tau=\operatorname{cone}\bigl(\tau_0\times\{0\},\Delta_{\tau,\rho}\times\{n_{\rho}\}\ (\rho\in p(\tau)(1))\bigr),\, \textrm{ where }\,
\Delta_{\tau,\rho}:=\{v\in N_{\Q}\mid(v,n_{\rho})\in\tau\}.
\]
\end{definition}

For such a cone $\tau$ and $\rho\in p(\tau)(1)$, the slice $\Delta_{\tau,\rho}$ is a non-empty polyhedron (the intersection of $\tau$ with an affine subspace) with tail $\tau_0$: if $v\in\Delta_{\tau,\rho}$ and $(w,0)\in\tau_0$ then $v+w\in\Delta_{\tau,\rho}$. Conversely, if $v+tw\in\Delta_{\tau,\rho}$ for all $t\geq0$, then $(w,0)=\lim_{t\to\infty}t^{-1}(v+tw,n_{\rho})\in\tau$. Moreover, $\tau_0$ is pointed. The Cayley fan is Cayley over $\Gamma_B$. For every subdivision $\Gamma$ of $\Gamma_B$,
\[
\Sigma(\D,B)\cap p^{-1}(\Gamma):=\{c\cap p^{-1}(\gamma)\mid c\in\Sigma(\D,B),\ \gamma\in\Gamma\}
\]
is a subdivision of $(\D,B)$ which is Cayley over $\Gamma$. Indeed, for cones $c,c'$ and $m=m_{1}+m_{2}\in(c\cap c')^{\vee}=c^{\vee}+c'^{\vee}$ ($m_{1}\in c^{\vee}$, $m_{2}\in c'^{\vee}$) one has $(c\cap c')\cap m^{\perp}=(c\cap m_{1}^{\perp})\cap(c'\cap m_{2}^{\perp})$. Applied to $c'=p^{-1}(\gamma)=N_{\Q}\times\gamma$, whose faces are the $N_{\Q}\times\gamma''$, $\gamma''\preceq\gamma$, this shows that the faces of $c\cap p^{-1}(\gamma)$ are of the form $c''\cap p^{-1}(\gamma'')$, and
\[
(c\cap p^{-1}(\gamma))\cap(c'\cap p^{-1}(\gamma'))=(c\cap c')\cap p^{-1}(\gamma\cap\gamma')
\]
is then a common face. Thus these cones form a fan refining $\Sigma(\D,B)$. For $c=\sigma_J$ and $\gamma\subset\Q^{J''}_{\geq0}$ in $\Gamma$, $p(\sigma_J\cap p^{-1}(\gamma))=\gamma\cap\Q^{J\cap J''}_{\geq0}$ is the face of $\gamma$ cut out by $\sum_{i\in J''\smallsetminus J}e_i$, hence lies in $\Gamma$. For $\lambda=\sum_{\rho}t_{\rho}n_{\rho}$ ($t_{\rho}\geq0$) in it, the slice description gives
\[
\sum_j\lambda_j\Delta_j+\omega=\sum_{\rho}t_{\rho}\left(\sum_j n_{\rho,j}\Delta_j\right)+\omega,
\]
so $\sigma_J\cap p^{-1}(\gamma)=\operatorname{cone}(\omega\times\{0\},(\sum_jn_{\rho,j}\Delta_j+\omega)\times\{n_{\rho}\})$ is generated by its tail and its slices at the $n_{\rho}$, i.e. is Cayley. Faces of $\sigma_J$, being Cayley cones of faces of $\omega$ and of the $\Delta_j$, are treated in the same way.

\subsection{The toroidal modification and its polyhedral divisors.}
A regular subdivision of the Cayley fan attached to the data of \cref{proposition: snc pp div and toroidal embeddings} gives a $T$-equivariant resolution of $\tilde{X}$ which, under a compatibility condition with the base, is covered by charts $\tilde{X}(\D_{\tau})$ for explicit polyhedral divisors $\D_{\tau}$ over a modification of $Y$. In what follows we describe these polyhedral divisors.

Let $\Sigma$ be a subdivision of $\Sigma(\D,B)$. For $\tau\in\Sigma$, put $Y_{\tau}:=Y\smallsetminus\bigcup_{i\notin J(\tau)}D_i$ and consider the presheaf of $M$-graded $\mathscr{O}_{Y_{\tau}}$-subalgebras of the constant sheaf $k(Y)[M]$ on $Y_{\tau}$
\begin{equation}\label{equation: cone algebra}
\mathscr{P}_{\tau}:U\longmapsto\bigoplus_{m\in M}\left(\sum_{\substack{a\in\Z^I \\ \ (m,a)\in\tau^{\vee}}}\H^0\left(U,\mathscr{O}_Y\left(-\sum_i a_iD_i\right)\right)\right)\chi^m.
\end{equation}
It is a subalgebra: $\tau^{\vee}$ is a semigroup, and the invertible sheaves $\mathscr{O}_Y(-\sum_ia_iD_i)$ multiply accordingly. Let $\mathscr{A}_{\tau}$ be the associated sheaf, and $\tilde{X}_{\tau}:=\Spec_{Y_{\tau}}(\mathscr{A}_{\tau})$ with projection $\pi_{\tau}$ and the $T$-action given by the grading. Note that \eqref{equation: cone algebra} depends only on $\tau$ and $B$, that only the $a_i$ with $i\in J(\tau)$ matter (since $D_i\cap Y_{\tau}=\emptyset$ otherwise), and that $\mathscr{A}(Y,\D)|_{Y_{\tau}}\subset\mathscr{A}_{\tau}$. Indeed, for $m\in\omega^{\vee}\cap M$, $\lfloor\D(m)\rfloor|_{Y_{\tau}}=-\sum_{i\in J(\tau)}a_iD_i$ with $a_i:=-\lfloor h_{\Delta_i}(m)\rfloor\geq-h_{\Delta_i}(m)$, and $(m,a)$ (with $a_i:=0$ for $i\notin J(\tau)$) lies in $\sigma_{J(\tau)}^{\vee}\subset\tau^{\vee}$ by Step 1 of the proof of \cref{proposition: snc pp div and toroidal embeddings}.

\begin{lemma}\label{lemma: local subdivision algebra}
Let $y\in Y$ be a closed point, $J:=\{i\mid y\in D_i\}$, and $V,y_j$ ($j\in J$), $\varphi$ as in Step 1 of the proof of \cref{proposition: snc pp div and toroidal embeddings} (these depend only on $B$). For $\tau\in\Sigma$, the cone $\tau_V:=\tau\cap(N_{\Q}\times\Q^J)$ is the face of $\tau$ cut out by $u:=(0,\sum_{i\notin J}e_i)$, and $\tau_V\in\Sigma_J:=\{\tau'\in\Sigma\mid\tau'\subset\sigma_J\}$. For every affine open $V'\subset V\cap Y_{\tau}$, the homomorphism of $M$-graded $\mathscr{O}(V')$-algebras
\[
\begin{split}
\theta:\mathscr{O}(V')\otimes_{k[z_j:j\in J]}k[\tau_V^{\vee}\cap(M\times\Z^J)]&\longrightarrow k(Y)[M],\\
f\otimes z^b\chi^m&\longmapsto fy^b\chi^m,
\end{split}
\]
is injective with image $\mathscr{P}_{\tau}(V')=\mathscr{A}_{\tau}(V')$. Consequently, $\mathscr{A}_{\tau}$ is quasi-coherent of finite type and, by \eqref{equation: local model}, $\tilde{X}_{\tau}\times_Y V\cong\pi^{-1}(V)\times_{X_{\sigma_J}}X_{\tau_V}$, $T$-equivariantly over $\pi^{-1}(V)$.
\end{lemma}
\begin{proof}
Since $\tau\subset\sigma_{J(\tau)}$, we obtain $\tau_V\subset\sigma_{J(\tau)}\cap(N_{\Q}\times\Q^J)=\sigma_{J\cap J(\tau)}\subset\sigma_J$. The morphism $(y_j)_j:V'\to\A^J$ is flat (it is $\varphi$ followed by a projection), so $\theta$, which is the base change of the inclusion $k[\tau_V^{\vee}\cap(M\times\Z^J)]\subset k[z_j^{\pm1}][M]$ followed by $\mathscr{O}(V')[y_j^{-1}][M]\subset k(Y)[M]$, is injective. On $V'$ we have $D_i|_{V'}=0$ for $i\notin J$ and $D_j|_{V'}=\div(y_j)|_{V'}$, so $\H^0(V',\mathscr{O}_Y(-\sum_ia_iD_i))=\mathscr{O}(V')y^{a_J}$ with $a_J:=(a_j)_{j\in J}$. Hence the degree-$m$ parts of $\mathscr{P}_{\tau}(V')$ and of the image of $\theta$ are $\sum_{a_J\in P}\mathscr{O}(V')y^{a_J}\chi^m$ and $\sum_{b\in Q}\mathscr{O}(V')y^b\chi^m$, where $P:=\{a_J\mid(m,a)\in\tau^{\vee}\}$ and $Q:=\{b\in\Z^J\mid(m,b)\in\tau_V^{\vee}\}$. As $\lambda_i\equiv0$ on $\tau_V$ for $i\notin J$, $P\subset Q$. Conversely, $\tau_V^{\vee}=\tau^{\vee}+\Q_{\geq0}(-u)=\tau^{\vee}+\Q u$ \cite[Section 1.2]{Ful93}, so for $b\in Q$ we can write $(m,b,0)=w-cu$ with $w\in\tau^{\vee}$, $c\in\Q$, and for an integer $l\geq\max(c,0)$ the vector $(m,b,0)+lu\in\tau^{\vee}$ is the lattice vector $(m,a)$ with $a_J=b$, and therefore $Q\subset P$. Finally, $V_{\tau}:=V\cap Y_{\tau}$ is affine (a principal open subset of $V$), and for affine $V''\subset V'\subset V_{\tau}$ the formula just proved gives $\mathscr{P}_{\tau}(V'')=\mathscr{O}(V'')\otimes_{\mathscr{O}(V')}\mathscr{P}_{\tau}(V')$ inside $k(Y)[M]$. Thus $\mathscr{P}_{\tau}$ agrees on affine open subsets of $V_{\tau}$ with the quasi-coherent sheaf associated with the finitely generated $\mathscr{O}(V_{\tau})$-algebra $\mathscr{P}_{\tau}(V_{\tau})$, and so does $\mathscr{A}_{\tau}$. The last assertion follows from \eqref{equation: local model} by associativity of the tensor product, since the spectrum of $\mathscr{O}(V)\otimes_{k[z_j]}k[\tau_V^{\vee}\cap(M\times\Z^J)]$ lies over $V_{\tau}$ ($y_i$ is invertible on it for $i\in J\smallsetminus J(\tau)$, since $(0,-e_i)\in\tau_V^{\vee}$).
\end{proof}

For $\tau'\preceq\tau$ in $\Sigma$, we have $\tau^{\vee}\subset\tau'^{\vee}$ and $J(\tau')\subset J(\tau)$, hence $\mathscr{A}_{\tau}\subset\mathscr{A}_{\tau'}|_{Y_{\tau}}$. Moreover, $\pi_{\tau'}(\tilde{X}_{\tau'})\subset Y_{\tau}$, since for $y\notin Y_{\tau}$ and $V$ as above there is $i\in J\smallsetminus J(\tau)$, so $y_i$ is invertible on $\tilde{X}_{\tau'}\times_Y V$ by the lemma. The induced $T$-equivariant morphism $\tilde{\iota}_{\tau'\tau}:\tilde{X}_{\tau'}\to\tilde{X}_{\tau}$ is, over each $V$, the base change of the open immersion $X_{\tau'_V}\subset X_{\tau_V}$ ($\tau'_V$ is a face of $\tau_V$), hence an open immersion. For $c,c'\preceq\tau$ one has $\tilde{X}_c\cap\tilde{X}_{c'}=\tilde{X}_{c\cap c'}$ inside $\tilde{X}_{\tau}$ (over $V$ this is $X_{c_V}\cap X_{c'_V}=X_{(c\cap c')_V}$). All these maps are induced by inclusions of subsheaves of $k(Y)[M]$, so the $\tilde{X}_{\tau}$, $\tau\in\Sigma$, glue along the $\tilde{X}_{\tau\cap\tau'}$ to a scheme $\tilde{X}_{\Sigma}$ with a $T$-action, a $T$-invariant morphism $\pi_{\Sigma}:\tilde{X}_{\Sigma}\to Y$ and, by $\mathscr{A}(Y,\D)|_{Y_{\tau}}\subset\mathscr{A}_{\tau}$, a $T$-equivariant morphism $\Phi_{\Sigma}:\tilde{X}_{\Sigma}\to\tilde{X}$ over $Y$. For $N=0$, a subdivision $\Gamma$ of $\Gamma_B$ and $\gamma\in\Gamma$, we write $Y_{\Gamma}:=\tilde{X}_{\Gamma}$, $Y_{\gamma}:=\tilde{X}_{\gamma}$, $\psi_{\Gamma}:=\Phi_{\Gamma}:Y_{\Gamma}\to Y$, and for $\rho\in\Gamma(1)$ let $E_{\rho}:=Y_{\Gamma}\smallsetminus\bigcup_{\gamma\in\Gamma,\ \rho\not\preceq\gamma}Y_{\gamma}$, with its reduced structure.

\begin{theorem}\label{theorem: combinatorial resolution}
Let $\Sigma$ be a subdivision of $\Sigma(\D,B)$.
\begin{enumerate}
\item \label{theorem: combinatorial resolution part i} In the situation of \cref{lemma: local subdivision algebra}, $\tilde{X}_{\Sigma}\times_Y V\cong\pi^{-1}(V)\times_{X_{\sigma_J}}X_{\Sigma_J}$ over $\pi^{-1}(V)$, where $X_{\Sigma_J}\to X_{\sigma_J}$ is the proper birational toric morphism of the refinement $\Sigma_J$ of $\sigma_J$.
\item \label{theorem: combinatorial resolution part ii} $\tilde{X}_{\Sigma}$ is a normal $T$-variety, $\Phi_{\Sigma}$ is $T$-equivariant, proper and birational, and an isomorphism over $U_B(\D)$. If $\Sigma$ is regular, then $\tilde{X}_{\Sigma}$ is smooth, so that $r'\circ\Phi_{\Sigma}$ is a $T$-equivariant resolution of singularities of $X(\D)$ in the situation of \cref{proposition: proper pullback modification} (regular subdivisions exist \cite[Section 2.6]{Ful93}).
\item \label{theorem: combinatorial resolution part iii} Let $\Gamma$ be a regular subdivision of $\Gamma_B$ and assume $\Sigma$ is Cayley over $\Gamma$. Then $Y_{\Gamma}$ is smooth, $\psi_{\Gamma}$ is proper and birational and an isomorphism over $Y\smallsetminus B$, each $E_{\rho}$ is smooth of pure codimension one with $\psi_{\Gamma}^{-1}(B)=\bigcup_{\rho}E_{\rho}$, the divisor $B_{\Gamma}:=\sum_{\rho}\sum_{F\subset E_{\rho}}F$ (sum over the irreducible components $F$ of the $E_{\rho}$) is a simple normal crossing divisor, and $Y_{\gamma}=Y_{\Gamma}\smallsetminus\bigcup_{\rho\notin\gamma(1)}E_{\rho}$. For $\tau\in\Sigma$ and $\gamma:=p(\tau)$ let $\D_{\tau}$ be the polyhedral divisor on $Y_{\gamma}$ with tail $\tau_0$,
\[
\D_{\tau}:=\sum_{\rho\in\gamma(1)}\sum_{F\subset E_{\rho}}\Delta_{\tau,\rho}\otimes(F\cap Y_{\gamma}),
\]
so that $\Supp(\D_{\tau})\leq B_{\Gamma}|_{Y_{\gamma}}$. Then the inclusions $\mathscr{A}_{\gamma}=\mathscr{A}_{\tau,0}\subset\mathscr{A}_{\tau}$ induce an affine $T$-invariant morphism $\tilde{X}_{\Sigma}\to Y_{\Gamma}$ over $Y$, and $\tilde{X}_{\tau}\cong\tilde{X}(\D_{\tau})$ $T$-equivariantly over $Y_{\gamma}$, compatibly with faces. Consequently, $W:=\tilde{X}_{\Sigma}$ is covered by the $T$-stable open subsets $\tilde{X}(\D_{\tau})$, with $\tilde{X}(\D_{\tau})\cap\tilde{X}(\D_{\tau'})=\tilde{X}(\D_{\tau\cap\tau'})$, where $(\tau\cap\tau')_0=\tau_0\cap\tau'_0$ and $\Delta_{\tau\cap\tau',\rho}=\Delta_{\tau,\rho}\cap\Delta_{\tau',\rho}$ (empty if $\rho\notin p(\tau\cap\tau')(1)$).
\end{enumerate}
\end{theorem}

\begin{proof}
(i) The cones of $\Sigma$ contained in $\sigma_J$ form a fan $\Sigma_J$ with support $\sigma_J$ (every $\tau\in\Sigma$ lies in some cone $c$ of $\Sigma(\D,B)$ and $\tau\cap\sigma_J=\tau\cap(c\cap\sigma_J)$ is a face of $\tau$), and $X_{\Sigma_J}\to X_{\sigma_J}$ is proper \cite[Theorem 3.4.11]{CLS11} and birational. By \cref{lemma: local subdivision algebra}, $\tilde{X}_{\tau}\times_Y V=\pi^{-1}(V)\times_{X_{\sigma_J}}X_{\tau_V}=\tilde{X}_{\tau_V}\times_Y V$ with $\tau_V\in\Sigma_J$, so both sides of (i) are unions of the open subschemes $\pi^{-1}(V)\times_{X_{\sigma_J}}X_{\tau'}$, $\tau'\in\Sigma_J$, glued along the same intersections.

(ii) Put $W_V:=\tilde{X}_{\Sigma}\times_Y V$. It is smooth over $X_{\Sigma_J}$ (a base change of the smooth morphism $p$), hence normal, with pairwise disjoint irreducible components, each mapping onto an open subset of $X_{\Sigma_J}$ and hence meeting the preimage of the torus $T_{N\times\Z^J}$. This preimage is $p^{-1}(T_{N\times\Z^J})=U_B(\D)\cap\pi^{-1}(V)$ (\cref{proposition: snc pp div and toroidal embeddings}\ref{proposition: snc pp div and toroidal embeddings part ii}), over which $X_{\Sigma_J}\to X_{\sigma_J}$ is an isomorphism. Thus $W_V$ is a normal variety and $\Phi_{\Sigma}$ is an isomorphism over $U_B(\D)\cap\pi^{-1}(V)$. Moreover, $W_V\to\pi^{-1}(V)$ is proper, being a base change of $X_{\Sigma_J}\to X_{\sigma_J}$. As the $\pi^{-1}(V)$ cover $\tilde{X}$, $\Phi_{\Sigma}$ is proper, so $\tilde{X}_{\Sigma}$ is separated and of finite type. It is irreducible, being a union of the irreducible open subsets $W_V$, any two of which meet in $U_B(\D)\cap\pi^{-1}(V\cap V')\neq\emptyset$ and the action is effective on $U_B(\D)$. If $\Sigma$ is regular, every cone of $\Sigma_J$ is generated by part of a basis of $\tilde{N}$ contained in the direct summand $N\times\Z^J$, hence by part of a basis of $N\times\Z^J$, so $X_{\Sigma_J}$ is smooth \cite[Theorem 3.1.19]{CLS11} and so is $W_V$. The last assertions follow from \cref{proposition: snc pp div and toroidal embeddings,proposition: proper pullback modification}.

(iii) Apply (i) and (ii) with $N=0$: $Y_{\Gamma}$ is smooth, $\psi_{\Gamma}$ is proper birational and an isomorphism over $Y\smallsetminus B=U_B(\D_0)$, and $Y_{\Gamma}\times_Y V=V\times_{\A^J}X_{\Gamma_J}$ with $g:V\times_{\A^J}X_{\Gamma_J}\to X_{\Gamma_J}\times\Gm^{n-r}$ étale (a base change of $\varphi$). By the orbit--cone correspondence, $E_{\rho}\times_Y V=V\times_{\A^J}D_{\rho}$ if $\rho\in\Gamma_J$ (since $D_{\rho}=X_{\Gamma_J}\smallsetminus\bigcup_{\gamma\not\succeq\rho}X_{\gamma}$) and $E_{\rho}\times_Y V=\emptyset$ otherwise, scheme-theoretically since $V\times_{\A^J}D_{\rho}$ is smooth over the reduced $D_{\rho}$. The preimage of $B\cap V=V(y_1\cdots y_r)$ in $X_{\Gamma_J}$ is $\bigcup_{\rho\in\Gamma_J(1)}D_{\rho}$ (each $n_{\rho}\neq0$ has a positive coordinate, and no $z_j$ vanishes on the torus), whence $\psi_{\Gamma}^{-1}(B)=\bigcup_{\rho}E_{\rho}$. For a cone $\gamma\in\Gamma_J$ of dimension $d$, its primitive generators $n_{\rho}$ ($\rho\in\gamma(1)$) extend to a basis of $\Z^J$ with dual basis $(u_{\rho})_{\rho\in\gamma(1)},(u'_l)_l$, so $k[\gamma^{\vee}\cap\Z^J]=k[z'_{\rho},z_l'^{\pm1}]$ with $z'_{\rho}:=\chi^{u_{\rho}}$, $X_{\gamma}\cong\A^d\times\Gm^{|J|-d}$, $D_{\rho}\cap X_{\gamma}=V(z'_{\rho})$ for $\rho\in\gamma(1)$ and $D_{\rho}\cap X_{\gamma}=\emptyset$ otherwise. Hence, over $V\times_{\A^J}X_{\gamma}$, each $E_{\rho}$ is $g^{-1}(V(z'_{\rho})\times\Gm^{n-r})$, smooth of pure codimension one, and the local equations $y'_{\rho}:=z'_{\rho}\circ g$ of the components of the $E_{\rho}$ through a point form part of a regular system of parameters ($g$ is étale). Thus $B_{\Gamma}$ is a simple normal crossing divisor, and $Y_{\gamma}=Y_{\Gamma}\smallsetminus\bigcup_{\rho\notin\gamma(1)}E_{\rho}$ since $X_{\gamma}=X_{\Gamma_J}\smallsetminus\bigcup_{\rho\notin\gamma(1)}D_{\rho}$.

Now let $\tau\in\Sigma$ and put $\gamma=p(\tau)$. Then $J(\gamma)=J(\tau)$ and $\gamma^{\vee}=\{a\mid(0,a)\in\tau^{\vee}\}$, so \eqref{equation: cone algebra} for $(N=0,\gamma)$ is the degree-zero part of \eqref{equation: cone algebra} for $\tau$, i.e. $\mathscr{A}_{\gamma}=\mathscr{A}_{\tau,0}$. This gives $T$-invariant morphisms $\tilde{X}_{\tau}\to Y_{\gamma}$ over $Y$, affine since $\tilde{X}_{\tau}$ and $Y_{\gamma}$ are affine over $Y_{\tau}$ and $Y_{\gamma}\to Y_{\tau}$ is separated. For $\tau'\preceq\tau$, the Cayley condition gives $p(\tau')\preceq p(\tau)$ ($p(\tau')$ is generated by the $n_{\rho}$ with $\Delta_{\tau',\rho}\neq\emptyset$, and $\Delta_{\tau',\rho}\subset\Delta_{\tau,\rho}$), and the morphisms are compatible with the chart inclusions, so they glue. Fix a closed point of $Y_{\gamma}$ over $y\in Y$, with $V,J$ as in \cref{lemma: local subdivision algebra}, so that $V'':=Y_{\gamma}\times_Y V=V\times_{\A^J}X_{\gamma_V}$ with $\gamma_V:=\gamma\cap\Q^J=p(\tau_V)$, and
\[
\H^0(\tilde{X}_{\tau}\times_Y V,\mathscr{O})=\mathscr{O}(V)\otimes_{k[z_j]}k[\tau_V^{\vee}\cap(M\times\Z^J)].
\]
The face $\tau_V$ of $\tau$, cut out by $(0,\sum_{i\notin J}e_i)$, which vanishes on $(v,n_{\rho})$ if and only if $n_{\rho}\in\Q^J$, is again Cayley over $\gamma_V$ with $\Delta_{\tau_V,\rho}=\Delta_{\tau,\rho}$ for $\rho\in\gamma_V(1)$ (a face of a Cayley cone cut out by a functional is the Cayley cone of the faces of its tail and slices, as in \cref{definition: Cayley subdivision}). With the coordinates $z'_{\rho},z_l'^{\pm1}$ on $X_{\gamma_V}$ as above, a vector $(m,b)\in M\times\Z^J$ lies in $\tau_V^{\vee}$ if and only if $m\in\tau_0^{\vee}$ and $\langle b,n_{\rho}\rangle=b_{\rho}\geq-h_{\Delta_{\tau,\rho}}(m)$ for all $\rho\in\gamma_V(1)$, where $b=\sum b_{\rho}u_{\rho}+\sum b_lu'_l$. Hence the degree-$m$ part of $k[\tau_V^{\vee}\cap(M\times\Z^J)]$ is zero for $m\notin\tau_0^{\vee}$ and, for $m\in\tau_0^{\vee}\cap M$, is the free $k[\gamma_V^{\vee}\cap\Z^J]$-module generated by $\prod_{\rho}z_{\rho}'^{-\lfloor h_{\Delta_{\tau,\rho}}(m)\rfloor}\chi^m$. Tensoring with $\mathscr{O}(V)$ over $k[z_j]$ gives
\[
\H^0(\tilde{X}_{\tau}\times_Y V,\mathscr{O})=
\bigoplus_{m\in\tau_0^{\vee}\cap M}\mathscr{O}(V'')\prod_{\rho\in\gamma_V(1)}y_{\rho}'^{-\lfloor h_{\Delta_{\tau,\rho}}(m)\rfloor}\chi^m\subset k(Y)[M].
\]
On the other hand, $V''$ is smooth, $\div(y'_{\rho})|_{V''}=E_{\rho}|_{V''}$, and the $E_{\rho}$ with $\rho\in\gamma(1)\smallsetminus\gamma_V(1)$ do not meet $V''$. Thus $\lfloor\D_{\tau}(m)\rfloor|_{V''}=\div(\prod_{\rho}y_{\rho}'^{\lfloor h_{\Delta_{\tau,\rho}}(m)\rfloor})|_{V''}$ and $\H^0(V'',\mathscr{A}(Y_{\gamma},\D_{\tau}))$ equals the right-hand side above. Since both $\tilde{X}_{\tau}$ and $\tilde{X}(\D_{\tau})$ are affine over $Y_{\gamma}$, and the $V''$ cover $Y_{\gamma}$, this gives $\tilde{X}_{\tau}\cong\tilde{X}(\D_{\tau})$ over $Y_{\gamma}$, equivariantly and compatibly with faces (all maps being inclusions of subsheaves of $k(Y)[M]$). Finally, $(\tau\cap\tau')\cap(N_{\Q}\times\{\lambda\})$ is the intersection of the slices of $\tau$ and $\tau'$ at $\lambda$, for $\lambda=0$ and $\lambda=n_{\rho}$. If $\rho\in p(\tau)(1)\cap p(\tau')(1)$ but $\rho\notin p(\tau\cap\tau')(1)$, then $n_{\rho}\notin p(\tau\cap\tau')$ (a ray of $\Gamma$ contained in a cone of $\Gamma$ is a ray of it), so the slice of $\tau\cap\tau'$ at $n_{\rho}$ is empty.
\end{proof}

We call the family $\S_W:=\{\D_{\tau}\mid\tau\in\Sigma\}$ of \cref{theorem: combinatorial resolution}\ref{theorem: combinatorial resolution part iii} the \emph{polyhedral divisors of the resolution}: the base is the toroidal modification $Y_{\Gamma}\to Y$ of $(Y,B)$ determined by $\Gamma$, whose boundary divisors $E_{\rho}$ are the strict transforms of the $D_i$ (for $\rho=\Q_{\geq0}e_i$) and exceptional divisors over the strata $D_{J(\rho)}$ (otherwise). The tail of $\D_{\tau}$ is $\tau_0$, its coefficient along every component of $E_{\rho}$ is the slice of $\tau$ at $n_{\rho}$, and the coefficients along the $E_{\rho}$ with $\rho\not\preceq p(\tau)$ are empty, so that $\Loc(\D_{\tau})=Y_{p(\tau)}$ in the sense of \cref{Section: resolution general case}. The family is closed under intersections and $W$ is glued from the $\tilde{X}(\D_{\tau})$ along the $\tilde{X}(\D_{\tau\cap\tau'})$, so $W$ plays the role of ``$\tilde{X}(\S_W)$''. This family is a divisorial fan in the sense of \cref{Section: resolution general case} only if each $\D_{\tau}$ is a pp-divisor and $Y_{p(\tau)}$ is the complement of the support of a semiample divisor, and even then $X(\D_{\tau})$ differs from $\tilde{X}(\D_{\tau})$ when $Y_{p(\tau)}$ is not affine (\cref{remark: nonaffine base}).

\begin{remark}\label{remark: affine base of resolution}
When $Y_{\Gamma}$ is affine, for instance in complexity one over an affine curve $Y$ (then $B$ is a finite set of points, $\Gamma_B=\{0\}\cup\{\Q_{\geq0}e_i\}$ has no proper subdivision, $Y_{\Gamma}=Y$, and every subdivision of $(\D,B)$ is Cayley over $\Gamma_B$, since a cone in $N_{\Q}\times\Q_{\geq0}e_i$ is generated by its slices at $0$ and $e_i$), these issues disappear: every effective divisor on $Y_{\Gamma}$ is Cartier and semiample with affine complement, so $Y_{p(\tau)}$ is an affine locus on which every $\Q$-divisor is semiample and big. Hence each $\D_{\tau}$ is a pp-divisor with $X(\D_{\tau})=\tilde{X}(\D_{\tau})$, $\S_W$ is a divisorial fan on $Y_{\Gamma}$ and $X(\S_W)\cong W$ equivariantly. Over a projective curve, the same holds after removing two points from the loci (every non-zero effective divisor being ample).
\end{remark}

\begin{remark}\label{remark: Cayley existence}
We do not know whether, for arbitrary $(\D,B)$, there exist a regular subdivision $\Gamma$ of $\Gamma_B$ and a regular subdivision of $\Sigma(\D,B)$ which is Cayley over $\Gamma$. The subdivision $\Sigma(\D,B)\cap p^{-1}(\Gamma)$ is Cayley but in general not regular, and refining it further (for instance, by a star subdivision at a point $(v,a)$ with $a$ interior to a cone of $\Gamma$ of dimension $\geq2$) may destroy the Cayley property. In general, modifications of the base are unavoidable. Let $N=\Z$, $Y=\A^2$, $B=D_1+D_2$ be the coordinate axes and $\Delta_1=\Delta_2=[\frac12,\infty)$. Then
\[
\sigma:=\sigma_{\{1,2\}}=\operatorname{cone}((1,0,0),(1,2,0),(1,0,2))\subset\Q\times\Q^2
\]
has multiplicity $4$, and no regular subdivision of $\Sigma(\D,B)$ is Cayley over $\Gamma_B$. Indeed, let $\Sigma$ be such a subdivision and $F:=\operatorname{cone}((1,2,0),(1,0,2))$, the face of $\sigma$ cut out by $(2,-1,-1)$. Some three-dimensional cone $\tau\in\Sigma$ has $\tau\cap F$ two-dimensional, and $p(\tau)=\Q^2_{\geq0}$. A point $(v,\lambda_1,\lambda_2)\in F$ with $\lambda_1,\lambda_2>0$ has $v=(\lambda_1+\lambda_2)/2$; writing it as $(w+\lambda_1v_1+\lambda_2v_2,\lambda_1,\lambda_2)$ with $w\in\tau_0\subset\Q_{\geq0}$ and $v_i\in\Delta_{\tau,e_i}\subset[\frac12,\infty)$ forces $w=0$ and $v_1=v_2=\frac12$. Hence $(1,2,0),(1,0,2)\in\tau$, and, being extremal rays of $\sigma\supset\tau$, they are rays of $\tau$; for any third primitive generator $(a,b,c)$, the determinant $4a-2b-2c$ is even, so $\tau$ is not regular. On the other hand, regular subdivisions of $\Sigma(\D,B)$ which are Cayley over the blow-up of $\Q^2_{\geq0}$ along the ray $(1,1)$ do exist (compare \cref{example: complexity two resolution}). We also do not know whether \cref{theorem: combinatorial resolution}\ref{theorem: combinatorial resolution part iii} extends to non-Cayley subdivisions with another choice of base, nor whether the proper modifications $\Phi_{\Sigma}$ can be chosen projective.
\end{remark}

\begin{example}\label{example: complexity two resolution}
Let $N=\Z$, $T=\Gm$, and let $X\subset\A^4=\Spec\left( k[x,y,s,t]\right)$ be the hypersurface $t^{2}=sxy(x-y)$. This hypersurface is normal and allows the following $T$-action 
\[ \lambda\cdot(x,y,s,t)=(x,y,\lambda^2s,\lambda t).\]
 Set $f:=xy(x-y)$. Then $X$ encoded by the pp-divisor $\D:=[\frac12,\infty)\otimes(D_1+D_2+D_3)$ on $Y:=\A^2$, where $D_1=\{x=0\}$, $D_2=\{y=0\}$, $D_3=\{x=y\}$, with tail $\Q_{\geq0}$. For every $m\in \Z$, its evaluations $\D(m)=\frac{m}{2}\div(f)$ are principal, and
\[
A[Y,\D]=\bigoplus_{m\geq0}k[x,y]f^{-\lfloor m/2\rfloor}\chi^m
\]
is the free $k[x,y]$-module with basis $s^j,ts^j$ ($j\geq0$), $t:=\chi$, $s:=f^{-1}\chi^2$, i.e. $k[x,y,s,t]/(t^2-fs)$ with $\deg(s)=2$, $\deg(t)=1$. The variety $X$ has complexity two and is singular along the lines $\{s=t=0,f=0\}$ and $\{x=y=t=0\}$. Let $\psi:Y'\to Y$ be the blow-up of the origin, where $B:=D_1+D_2+D_3$ is not simple normal crossing, with exceptional curve $E$ and strict transforms $D'_i$; then $B':=\sum D'_i+E$ is simple normal crossing with $D'_i\cap E=\{P_i\}$, $P_1,P_2,P_3$ distinct. By \cref{lemma: birational modification of a pp-div}, $\D':=\psi^*\D$ has coefficients $[\frac12,\infty)$ along $D'_i$ and $[\frac32,\infty)$ along $E$ ($f$ vanishes to order $3$ at the origin). Writing vectors of $\Q\times\Q^{\{E,i\}}$ as $(v,\lambda_E,\lambda_i)$, the Cayley cones are
\[
\begin{gathered}
\sigma_{\{i\}}=\operatorname{cone}((1,0,0),(1,0,2)),\qquad
\sigma_{\{E\}}=\operatorname{cone}((1,0,0),(3,2,0)),\\
\sigma_{\{E,i\}}=\operatorname{cone}((1,0,0),(3,2,0),(1,0,2)),
\end{gathered}
\]
the last cone having multiplicity $4$. Let $\Gamma$ be the regular subdivision of $\Gamma_{B'}$ inserting the ray $\Q_{\geq0}(e_E+e_i)$ in each quadrant $\operatorname{cone}(e_E,e_i)$. Then $Y'':=Y'_{\Gamma}$ is the blow-up of $Y'$ at $P_1,P_2,P_3$ (by \cref{theorem: combinatorial resolution}\ref{theorem: combinatorial resolution part i}, blowing up commutes with flat base change), with boundary divisors $D''_i,E''$ (strict transforms) and the exceptional curves $F_i$ over $P_i$. For each $i$, subdivide $\sigma_{\{E,i\}}$ into the four cones
\[
\begin{aligned}
\tau^i_1&:=\operatorname{cone}((1,0,0),(2,1,0),(2,1,1)),\\
\tau^i_2&:=\operatorname{cone}((2,1,0),(3,2,0),(2,1,1)),\\
\tau^i_3&:=\operatorname{cone}((1,0,0),(1,0,1),(2,1,1)),\\
\tau^i_4&:=\operatorname{cone}((1,0,1),(1,0,2),(2,1,1)),
\end{aligned}
\]
of determinant $\pm1$. Here $(2,1,1)$ spans the slice $[2,\infty)$ of $\sigma_{\{E,i\}}$ at $e_E+e_i$, and $(2,1,0),(1,0,1)$ subdivide the faces $\sigma_{\{E\}},\sigma_{\{i\}}$ into regular cones, compatibly for $i=1,2,3$. Thus these cones and their faces form a regular subdivision $\Sigma$ of $(\D',B')$, which is Cayley over $\Gamma$ (each cone is generated by its tail and its slices, computed below). By \cref{theorem: combinatorial resolution}, $W:=\tilde{X}'_{\Sigma}\to\tilde{X}'\to X$ is a $\Gm$-equivariant resolution of singularities, $W\to Y''$ is affine, and $W$ is covered by the $\tilde{X}(\D_{\tau})$, $\tau\in\Sigma$ maximal, where, with
\[
U_{i,1}:=Y''\smallsetminus\left(\bigcup_jD''_j\cup\bigcup_{j\neq i}F_j\right),\qquad
U_{i,2}:=Y''\smallsetminus\left(E''\cup\bigcup_{j\neq i}(D''_j\cup F_j)\right),
\]
\[
\begin{array}{llll}
\D_{\tau^i_1}=[2,\infty)\otimes E''+[2,\infty)\otimes F_i,&\text{tail }\Q_{\geq0},&\text{on }U_{i,1},\\[3pt]
\D_{\tau^i_2}=[\frac32,2]\otimes E''+\{2\}\otimes F_i,&\text{tail }\{0\},&\text{on }U_{i,1},\\[3pt]
\D_{\tau^i_3}=[1,\infty)\otimes D''_i+[2,\infty)\otimes F_i,&\text{tail }\Q_{\geq0},&\text{on }U_{i,2},\\[3pt]
\D_{\tau^i_4}=[\frac12,1]\otimes D''_i+\{2\}\otimes F_i,&\text{tail }\{0\},&\text{on }U_{i,2},
\end{array}
\]
All other coefficients are equal to the tail.
\end{example}

\end{document}